\documentclass[a4paper, 12pt]{article}
\usepackage{amsmath,amsthm,amsfonts,amssymb,stmaryrd,amscd}
\usepackage[top=2cm, bottom=2.5cm, left=2.2cm, right=2.2cm]{geometry}
\usepackage[french, english]{babel}

\input xy
\xyoption{all}

\usepackage{color}
\usepackage{makeidx}
\usepackage{hyperref}
\makeindex 

\numberwithin{equation}{section}

{\theoremstyle{definition}
\newtheorem{definition}{Definition}[section]

\newtheorem{notation}[definition]{Notation}

\newtheorem{remark}[definition]{Remark}
\newtheorem{facts}[definition]{Facts}

}
\newtheorem{proposition}[definition]{Proposition}
\newtheorem{proposition-definition}[definition]{Proposition-Definition}
\newtheorem{lemma}[definition]{Lemma}
\newtheorem{theorem}[definition]{Theorem}
\newtheorem{corollary}[definition]{Corollary}

\newcommand{\R}{\mathbb{R}}

\newcommand{\cG}{\mathcal{G}}
\newcommand{\cJ}{\mathcal{J}}
\newcommand{\cE}{\mathcal{E}}

\newcommand{\cL}{\mathcal{L}}

\newcommand{\cK}{\mathcal{K}}
\newcommand{\C}{\mathbb{C}}

\newcommand{\ev}{{\rm ev}}
\newcommand{\dom}{\mathrm{dom}}

\newcommand{\cA}{\mathcal{A}}

\newcommand{\cP}{\mathcal{P}}

\newcommand{\cM}{\mathcal{M}}

\newcommand{\cS}{\mathcal{S}}

\newcommand{\resp}{{\it resp.}\/ }
\newcommand{\id}{{\hbox{id}}}

\newcommand{\ie}{{\it i.e.}\/ }

\newcommand{\cf}{{\it cf.}\/ }

\newcommand{\gA}{\mathfrak{A}}

\newcommand{\gH}{\mathfrak{H}}
\newcommand{\gag}{\mathfrak{g}}

\def\gpd{\,\lower1pt\hbox{$\longrightarrow$}\hskip-.24in\raise2pt
             \hbox{$\longrightarrow$}\,}

\newcommand{\To}{\longrightarrow}
\reversemarginpar

\everymath{\displaystyle}

\begin{document}

\begin{center}
{\Large \bf  Pseudodifferential extensions and adiabatic deformation of smooth groupoid actions}

\bigskip
{\sc by Claire Debord and Georges Skandalis}

\end{center}

{\footnotesize
Laboratoire de Math\'ematiques, UMR 6620 - CNRS
\vskip-4pt Universit\'e Blaise Pascal, Campus des C\'ezeaux, BP {\bf 80026}
\vskip-4pt F-63171 Aubi\`ere cedex, France
\vskip-4pt claire.debord@math.univ-bpclermont.fr

\vskip 2pt Universit\'e Paris Diderot, Sorbonne Paris Cit\'e
\vskip-4pt  Sorbonne Universit\'es, UPMC Paris 06, CNRS, IMJ-PRG
\vskip-4pt  UFR de Math\'ematiques, {\sc CP} {\bf 7012} - B\^atiment Sophie Germain 
\vskip-4pt  5 rue Thomas Mann, 75205 Paris CEDEX 13, France
\vskip-4pt skandalis@math.univ-paris-diderot.fr
}
\bigskip

\centerline{\bf Abstract}

The adiabatic groupoid $\cG_{ad}$ of a smooth groupoid $\cG$ is a deformation relating $\cG$ with its algebroid. In a previous work, we constructed a natural action of $\R$ on the C*-algebra of zero order pseudodifferential operators on $\cG$ and identified the crossed product with a natural ideal $J(\cG)$ of $C^*(\cG_{ad})$. In the present paper we show that $C^*(\cG_{ad})$ itself is a pseudodifferential extension of this crossed product in a sense introduced by Saad Baaj. Let us point out that we prove our results in a slightly more general situation:  the smooth groupoid $\cG$ is assumed to act on a C*-algebra $A$. We construct in this generalized setting the extension of order $0$ pseudodifferential operators $\Psi(A,\cG)$  of the associated crossed product $A\rtimes \cG$. We show that $\R$ acts naturally on $\Psi(A,\cG)$ and identify the crossed product of $A$ by the action of the adiabatic groupoid $\cG_{ad}$ with an extension of the crossed product $\Psi(A,\cG)\rtimes \R$.  Note that our construction of $\Psi(A,\cG)$ unifies the ones of Connes (case $A=\C $) and of Baaj ($\cG$ is a Lie group).


\renewcommand\theenumi{\alph{enumi}}
\renewcommand\labelenumi{\rm {\theenumi})}
\renewcommand\theenumii{\roman{enumii}}

\section{Introduction}

Alain Connes in \cite[{Chap. VIII}]{ConnesLNM} pointed out that smooth groupoids offer a perfect setting for index theory. Since then, this fact has been explored and exploited by Connes as well as many other authors, in many geometric situations (see \cite{DebordLescureIndex} for a review).

In \cite[{section II.5}]{ConnesNCG}, A. Connes constructed a beautiful groupoid, that he called the ``tangent groupoid'', which interpolates between the pair groupoid $M\times M$ of a (smooth, compact) manifold $M$ and the tangent bundle $TM$ of $M$. He showed that this groupoid describes the  analytic index on $M$  in a way not involving (pseudo)differential operators at all, and gave a proof of the Atiyah-Singer Index Theorem based on this groupoid. 

This idea of a deformation groupoid was then used in \cite[{section III}]{HilSkFeuill}, and extended in \cite{MonthPie, NWX} to the general case of a smooth groupoid, where the authors associated to every smooth groupoid $\cG$ an \emph{adiabatic groupoid} $\cG_{ad}$, which is obtained by applying the ``deformation to the normal cone'' construction to the inclusion $\cG^{(0)}\to \cG$ of the unit space of $\cG$ into $\cG$. Moreover, it was shown in \cite[{Th\'eor\`eme 2.1}]{MonthPie} that this adiabatic groupoid still describes the analytic index of the groupoid $\cG$ in this generalized situation.

In \cite{DS1}, we further explored the relationship between pseudodifferential calculus on $\cG$ and its adiabatic deformation $\cG_{ad}$. An ideal $J(\cG)\subset C^*(\cG_{ad})$ which sits in an exact sequence $0\to J(\cG)\to C^*(\cG_{ad})\to C(\cG^{(0)})\to 0$ plays a crucial role in our constructions. We construct a canonical Morita equivalence between the algebra $\Psi^*(\cG)$ of order $0$ pseudodifferential operators on $\cG$ and the crossed product $J(\cG)\rtimes \R_+^*$ of $J(\cG)$ by the natural action of $\R_+^*$. \\ It appeared that $J(\cG)$ is canonically isomorphic to the crossed product $\Psi^*(\cG)\rtimes \R$ associated with a natural action of $\R$ on the algebra $\Psi^*(\cG)$. A natural question is then: can one recognize the $C^*$-algebra $C^*(\cG_{ad})$ in these terms ? 

In the present paper, we answer this question, thanks to  \cite{BaajOPD1, BaajOPD2}, where Baaj constructed {an extension of pseudodifferential operators of order $0$ of the crossed product}  of a $C^*$-algebra $A$ by the action of a Lie group $H$ - with Lie algebra $\gH$. Denote by $S^*\gH$ the sphere in $\gH^*$.
Baaj's exact sequence reads $$0\to A\rtimes H\To \Psi_0^*(A,H)\overset{\sigma}{\To}C(S^*\gH)\otimes A\to 0.$$

Let $\mu :C(\cG^{(0)})\to \Psi^*(\cG)$ be the inclusion by multiplication operators. In the present paper, we construct a commutative diagram, whose first line is Baaj's exact sequence:
 $$\xymatrix{
0\ar[r]&\Psi^*(\cG)\rtimes \R\ar[r]  &\Psi_0 ^*(\Psi^*(\cG),\R)\ar[r]^{\sigma\ \ \ }& \Psi^*(\cG)\oplus \Psi^*(\cG)\ar[r]& 0\\
0\ar[r]& J(\cG)\ar[r]\ar[u]_{\simeq}&C^*(\cG_{ad})\ar[r]\ar[u] &C(\cG^{(0)})\ar[u]_{\mu _0}\ar[r]& 0 
}\eqno (1)
$$ 
where $\mu_0(f)=(\mu(f),0)$.

Moreover, we show that all the morphisms of the above diagram are equivariant with respect to the natural actions of $\R_+^*$: \begin{itemize}
\item We consider $\R_+^*$ as the dual group of $\R$ and thus it acts on the crossed product $\Psi^*(\cG)\rtimes \R$ via the dual action. This dual action extends (uniquely) to Baaj's pseudodifferential extension $\Psi ^*_0(\Psi^*(\cG),\R)$ and is trivial at the quotient level. 

\item The action of $\R_+^*$ on the second line is the canonical action  on the adiabatic groupoid  by the natural rescaling, and the crossed product $C^*(\cG_{ad})\rtimes \R_+^*$ is the $C^*$-algebra of the ``gauge adiabatic groupoid'' $\cG_{ga}$ considered in \cite{DS1}. 
\end{itemize}
In particular, this allows us to give also a description of the algebra $C^*(\cG_{ga})$ as a pseudodifferential extension.

As a side construction, we define the pseudodifferential extension of an action $\alpha$ of a smooth groupoid $\cG$ - in the setting introduced by Le Gall in  \cite{PYLG, PYLG2}. {This is a short exact sequence $$0\to A\rtimes _\alpha \cG \longrightarrow \Psi^*(A,\alpha,\cG)\overset{\sigma_\alpha}\longrightarrow A\otimes _{C_0(M)}C(S^* \gA \cG)\to 0.\eqno(2)$$} This construction generalizes both the pseudodifferential calculus on a smooth groupoid of \cite{ConnesLNM, ConnesNCG, MonthPie, NWX} and the pseudodifferential calculus of a crossed product by a Lie group of \cite{BaajOPD1, BaajOPD2}. Our main result, Theorem \ref{maintheorem}, is stated (and proved) in this general frame: in diagram (1) we allow the groupoid $\cG$ to act on a $C^*$-algebra $A$ and replace groupoid $C^*$-algebras by crossed products. We should note that the connecting map of extension (2) is the analytic index in this context. In the same way as in \cite{MonthPie, NWX}, the crossed product by the adiabatic groupoid allows to define the analytic index too.

\medskip 
Here are some examples of natural actions of smooth groupoids which are relevant to our constructions.
\begin{enumerate}\renewcommand\theenumi{\arabic{enumi}}
\renewcommand\labelenumi{{\theenumi}.}
\item Already an interesting case appears when $A=C_0(X)$ where $X$ is a smooth manifold, endowed with a smooth submersion $p:X\to M=\cG^{(0)}$ and $\cG$ acts on the fibers. The action of $\cG$ is given by a diffeomorphism $\alpha:\cG{}_s\times_pX\to X{}_p\times_r\cG$ of the form $(\gamma,x)\mapsto (\alpha_\gamma(x),\gamma)$, which satisfies $\alpha_{\gamma_1\gamma_2} =\alpha_{\gamma_1}\alpha_{\gamma_2}$. Here, $\cG{}_s\times_pX$ is a smooth groupoid $\cG_X$ with objects $X$, source and range maps given by $s(\gamma,x)=x$, $r(\gamma,x)=\alpha_\gamma(x)$ composition $(\gamma',\alpha_\gamma(x))(\gamma,x)=(\gamma'\gamma,x)$ and inverse $(\gamma,x)^{-1}=(\gamma^{-1},\alpha_\gamma(x))$. In that case, the crossed product $A\rtimes_\alpha \cG$, the extension $\Psi^*(A,\alpha,\cG)$, the crossed product $(A\otimes \R_+)\rtimes \cG_{ad}$ identify respectively with the groupoid $C^*$-algebra $C^*(\cG_X)$, the pseudodifferential extension $\Psi^*(\cG_X)$ and the $C^*$-algebra $C^*((\cG_X)_{ad})$ of the adiabatic deformation of the groupoid $\cG_X$.
\item Let $G$ be a Lie group acting on a $C^*$-algebra $A$. The corresponding adiabatic and gauge adiabatic deformations of $G$ are groupoids with objects $\R_+$. They naturally act on the $C_0(\R_+)$ algebra $A\otimes C_0(\R_+)$ - and the associated action is an important piece in our constructions - see section \ref{adgpdaction}.
\item An interesting family of examples of groupoid actions comes from $1$-cocycles (generalized morphisms in the sense of \cite[section I]{HilSkFeuill}, \cite[{Section 2.2}]{PYLG}) of a groupoid $\cG$ to a Lie group. For instance, an equivariant vector bundle is equivalent to a cocycle from $\cG$ to $GL_n(\R)$.  Then every algebra $A$ endowed with an action of $G$ gives rise to a $\cG$-algebra. This construction is studied in \cite{PYLG} where several examples connected with $K$-theory and index theory are studied.  The corresponding pseudodifferential extension and associated actions of the adiabatic groupoid appear very naturally in this context.
\end{enumerate}

\medskip The paper is organized as follows:

In the second section, we briefly review the action of a locally compact groupoid and the corresponding full and reduced crossed products (\cf \cite{PYLG, PYLG2, Tu1, Tu2, Pat}).

In the third section, we review Baaj's construction and discuss the dual action.

In the fourth section we generalize Baaj's construction  to the case of actions of smooth groupoids. 

The fifth section establishes the above mentioned equivariant commutative diagram.

Finally, we gathered a few rather well known facts on unbounded multipliers in an appendix.


\begin{notation}
If $A$ is a $C^*$-algebra, we denote by $\cM(A)$ its multiplier algebra. 

Recall that, if $A$ and $B$ are $C^*$-algebras, a morphism $f:A\to \cM(B)$ is said to be non degenerate if $f(A).B=B$; a non degenerate morphism extends uniquely to a morphism $\tilde f:\cM(A)\to \cM(B)$ - this extension is strictly continuous (\ie continuous with respect to the natural topologies of the multipliers).

Recall that an ideal $J$ of a $C^*$-algebra $A$ is said to be essential if the morphism $A\to \cM(J)$ is injective, \ie if $a\in A$ is such that $aJ=\{0\}$ then $a=0$.
\end{notation}

\begin{remark}\label{essential}
Note that if $\pi: A\to B$ is a surjective morphism of $C^*$-algebras and $J$ an essential ideal in $B$ then $\pi^{-1} (J)$ is essential in $A$.
\end{remark}

\section{Actions of locally compact groupoids and crossed products}

In this section we briefly recall a few facts about actions of locally compact groupoids and the corresponding crossed products as defined by Le Gall in \cite{PYLG, PYLG2}. See also \cite{Tu1, Tu2, Pat}.

\subsection{Actions of locally compact groupoids}\label{essonfrerK1}

\subsubsection{${\mathbf C_0(X)}$-algebras}\label{soeur}

\begin{description}
\item[${\mathbf C_0(X)}$-algebras.] Recall (\cite{DG}, \cite[{Def. 1.5}]{Kasparov}) that if $X$ is a locally compact space, a $C_0(X)$-algebra is a pair $(A,\theta)$, where $A$ is a $C^*$-algebra and $\theta$ is a non degenerate $*$-homomorphism  $\theta:C_0(X)\to Z\cM(A)$ from $C_0(X)$ to the center of the multiplier algebra of $A$. 

\item[Fibers.] If $A$ is a $C_0(X)$-algebra, we define its fiber $A_x$ for every point $x\in X$ by setting $A_x=A/C_xA$ where $C_x = \{h\in C_0 (X);\ h(x) = 0\}$. Let $a\in A$ and denote by $a_x\in A_x$ its class; we have $\|a\|=\sup_{x\in X} \|a_x\|$.  In particular $a $ is completely  determined by the family $(a_x)_{x\in X}$  and the bundle $A$ is semi-continuous in  the sense that for all $a\in A$ the map $x\mapsto \|a_x\|$ is upper semi-continuous. 

\item [${\mathbf C_0(X)}$-morphisms.] A $C_0(X)$-linear homomorphism $\alpha:A\to B$ of $C_0(X)$-algebras determines for each $x\in X$ a $*$-homomorphism $\alpha_x:A_x\to B_x$. Since $\alpha(a)$ is determined by the family $(\alpha(a))_x=\alpha_x(a_x)$, the morphism $\alpha$ is determined by the family $(\alpha_x)_{x\in X}$.

\item[Restriction to locally closed sets; pull back.] More generally, if $U\subset X$ is an open subset, we define the $C_0(U)$-algebra $A_U$ by putting $A_U=C_0(U)A$; if $F\subset X$ is a closed subset, we define the $C_0(F)$-algebra $A_F=A/A_{X\setminus F}$; if $Y=U\cap F$ is a locally closed subset of $X$ we  put $A_Y=(A_U)_Y$ (which is canonically isomorphic to $(A_F)_Y$).

Recall that if $f:Y\to X$ is a continuous map between locally compact spaces and $A$ is a $C_0(X)$-algebra, we may define $f^*(A)$ in the following way: we restrict the $C_0(X\times Y)$-algebra $A\otimes C_0(Y)$ to the graph $\{(x,y)\in X\times Y;\ f(y)=x\}$ of $f$ which is a closed subset of $X\times Y$ canonically homeomorphic with $Y$.
\end{description}

\begin{notation}
 As $f^*(A)$ is a quotient of $A\otimes C_0(Y)$, we have a non degenerate morphism  $a\mapsto a\circ f$ from $A$ to the multiplier algebra of $f^*(A)$, where $a\circ f$ is the image of $a\otimes 1$ in the quotient $f^*(A)$ of $A\otimes C_0(Y)$.
\end{notation}

\subsubsection{Actions of groupoids}

\begin{definition}(\cite[Definition 2.2]{PYLG2}). Let $\cG$ be a locally compact groupoid with basis $X$. A continuous action of $\cG$ on a $C_0(X)$-algebra $A$ is an isomorphism of $C_0(\cG)$-algebras $\alpha : s^*A \to r^*A$ such that, for all $(\gamma_1,\gamma_2)\in \cG^{(2)}$ we have $\alpha_{\gamma_1\gamma_2}=\alpha_{\gamma_1}\circ\alpha_{\gamma_2}$.
\end{definition}

\begin{remark}
An action of a non Hausdorff groupoid $\cG$ on a $C_0(X)$-algebra $A$ (with $X=\cG^{(0)}$) is given by isomorphisms $\alpha_U:s_U^*(A)\to r_U^*(A)$ for every Hausdorff open subset $U$ of $X$ - where $s_U,r_U$ are the restrictions of $r$ and $s$ to $U$. These isomorphisms must agree on the intersection $U\cap V$ of two such sets. It follows that the family $(r_U)$ gives rise to isomorphisms $\alpha_\gamma:A_{s(\gamma)}\to A_{r(\gamma)}$ for $\gamma\in \cG$. We further impose that these isomorphisms satisfy $\alpha_{\gamma_1\gamma_2}=\alpha_{\gamma_1}\circ\alpha_{\gamma_2}$ for all $(\gamma_1,\gamma_2)\in \cG^{(2)}$.\\
In the sequel of the paper, we will consider Hausdorff groupoids for simplicity of the exposition. Nevertheless, all our constructions and results extend in the usual way to the non Hausdorff case \cite[{section 6}]{ConnesSurvey}, see also \cite[{section I.B}]{KhoshSk}. Note that the non trivial part of any kind of pseudodifferential calculus concentrates in a Hausdorff neighborhood of the space of units. 
\end{remark}

\subsection{Crossed products}

The (full and reduced) crossed product $A\rtimes_\alpha \cG$ of an action $\alpha $ of a groupoid $\cG$ with (right) Haar system $(\nu ^x)_{x\in X}$ on a $C^*$-algebra $A$ is defined in \cite{PYLG2, Pat2}. Let us briefly recall these constructions.

\subsubsection{The full crossed product}
The vector space $C_c(r^*A)=C_c(\cG).r^*(A)$ of elements of $r^*A$ with compact support is naturally a convolution $*$-algebra. For $f,g\in  C_c(r^*A)$ and $\gamma\in \cG$, we have $$(f\ast g)_\gamma=\int_{\cG^{r(\gamma)}}f_{\gamma_1}\alpha_{\gamma_1}(g_{\gamma_1^{-1}\gamma})\,d\nu ^{r(\gamma)}(\gamma_1)\quad \hbox{and} \quad (f^*)_\gamma=\alpha_{\gamma}^{-1}(f_{{\gamma}^{-1}})$$

There is a $\|\ \|_1$ norm given by $$\|f\|_1=\sup_{x\in X} \max\left(\int_{\cG^x}\|f_\gamma\|d\nu ^x(\gamma),\int_{\cG^x}\|f_{\gamma^{-1}}\|d\nu ^x(\gamma)\right)$$on this algebra and the corresponding completion is a Banach $*$-algebra $L^1(r^*A,\nu )$ (recall that $X$ is the basis $\cG^{(0)}$ of $\cG$).

The full crossed product $A\rtimes_\alpha \cG$ is the enveloping $C^*$-algebra of $L^1(r^*A,\nu )$. The algebras $A$ and $C^*(\cG)$ sit in the multipliers of $A\rtimes_\alpha \cG$ in a non degenerate way, and $A\rtimes_\alpha \cG$ is the closed vector span of products $a.f$ with $a\in A$ and $f\in C^*(\cG)$. Note that $C_0(X)$ sits both in the multipliers of $C^*(\cG)$ and of $A$; its images in $\cM(A\rtimes_\alpha \cG)$ agree.  

\subsubsection[Covariant representations]{Covariant representations (see \cite[p. 1466]{Pat2} -  see also \cite[{section II.1}]{Ren})}
 The representations of $A\rtimes _\alpha \cG$ can easily be described as in \cite[{Theorem 1.21, p. 65}]{Ren}. Such a representation gives rise to representations of $A$ and $C^*(\cG)$. We thus obtain:
\begin{itemize}
\item The representation of $C_0(X)$  corresponds to a measure $\mu$ on $X$ and a measurable field of Hilbert spaces $(H_x)_{x\in X}$.
\item The representation of the $C_0(X)$-algebra $A$ is given by a measurable family $\pi=(\pi_x)_{x\in X}$ where $\pi_x:A_x\to \cL(H_x)$ is a $*$-representation.

\item The representation of $C^*(\cG)$ gives rise to a representation of $\cG$ in the sense of \cite[{def. 1.6, p. 52}]{Ren}. In other words, the measure $\mu$ is quasi-invariant (\ie $\mu\circ\nu $ is quasi-invariant by the map $\gamma\mapsto \gamma^{-1}$) and we have a measurable family $U=(U_\gamma)_{\gamma\in \cG}$ where $U_\gamma:H_{s(\gamma)}\to H_{r(\gamma)}$ is (almost everywhere) unitary and satisfies (almost everywhere) $U_{\gamma_1\gamma_2}=U_{\gamma_1}U_{\gamma_2}$.

\item The covariance property then reads: $\pi_{r(\gamma)}\circ \alpha_\gamma=Ad_{U_\gamma}\circ \pi_{s(\gamma)}$  (almost everywhere).
\end{itemize}

Conversely,  such data $(\mu,H,\pi,U)$ can be integrated to a representation of  $A\rtimes _\alpha \cG$.

\subsubsection[The reduced crossed product]{The reduced crossed product (see \cite{Ren, KhoshSk})}
The reduced crossed product $A\rtimes_{\alpha,red} \cG$ is the quotient of $A\rtimes_\alpha \cG$ corresponding to the family of regular representations on the Hilbert modules $A_x\otimes L^2(\cG^x;\nu ^x)$ for $x\in X$.

If $\cG$ is amenable (see \cite{ADR} {for a discussion on amenability of groupoids}) then the morphism $A\rtimes_\alpha \cG\to A\rtimes_{\alpha,red} \cG$ is an isomorphism.

The reduced crossed product has a faithful representation on the Hilbert $A$-module $\cE=L^2(\cG;\nu )\otimes _{C_0(X)}A$ where $L^2(\cG;\nu )$ is the Hilbert $C_0(X)$ module described in \cite[{Theorem 2.3}]{KhoshSk} (if $\cG$ is Hausdorff). The module $\cE$ is the completion of $C_c(\cG;s^*A)$ with respect to the $A$-valued inner product satisfying $(\langle \xi|\eta\rangle )_x=\int _{\cG_x}\xi ^*_\gamma\eta_\gamma d\nu _x(\gamma)$ \big(where $(\nu_x)_{x\in X}$ is the corresponding left Haar system given by $\int f(\gamma) d\nu_x(\gamma)=\int f(\gamma^{-1}) d\nu ^x(\gamma)$) and, right action given by $(\xi a)_\gamma=\xi_\gamma a_{s(\gamma)}$\big).  

Denote by $\lambda $ the action of  $C_{red}^*(\cG)$ by (left) convolution on the Hilbert $C_0(X)$-module $L^2(\cG;\nu )$; the left action of $C^*(\cG)$ is given by $f\mapsto \lambda(f)\otimes _{C_0(X)}1$. The action of  $A$ is given by $a.\xi=\Big(\alpha^{-1}(a\circ r)\Big)\xi$: in other terms $(a.\xi)_\gamma=\alpha_\gamma^{-1}(a_{ r(\gamma)})\xi_\gamma$.

It follows, that if $\pi =\int _{X}^\oplus \pi_x\,d\mu(x)$ is a faithful representation of $A$, the corresponding representation of $A\rtimes _{\alpha, red}\cG$ on $\int_X^\oplus L^2(\cG_x,\nu_x)\otimes H_x\,d\mu(x)$ is faithful.

\subsubsection[Invariant ideals and exact sequences]{Invariant ideals and exact sequences (see  \cite[Theorem 3]{Pat2})}

Let $J\subset A$ be an ideal in $A$. Note that both $J$ and $A/J$ are then $C_0(X)$ algebras - recall that $X=\cG^{(0)}$. Assume that $J$ is invariant under the action of $\cG$ which means that $\alpha(s^*(J))=r^*(J)$. Then $\alpha$ yields actions of $\cG$ on $J$ and $A/J$. 
\begin{lemma} \cite[Theorem 3]{Pat2}
We have an exact sequence of \emph{full} crossed products: $$0\to J\rtimes_\alpha \cG\to A\rtimes_\alpha \cG\to (A/J)\rtimes_\alpha \cG\to 0.$$
\begin{proof}
The only thing which is not completely obvious in this sequence is that the morphism $(A\rtimes_\alpha \cG)/(J\rtimes_\alpha \cG)\to A/J\rtimes_\alpha \cG$ is injective. To see that, take a faithful representation of $(A\rtimes_\alpha \cG)/(J\rtimes_\alpha \cG)$; it is a covariant representation of $A$ and $\cG$ which vanishes on $J$, and therefore a covariant representation of $A/J$ and $\cG$.
\end{proof}
\end{lemma}

If $J$ is a $\cG$-invariant essential ideal in $A$, then at the level of \emph{reduced} crossed products, the ideal $J\rtimes_{\alpha,red}\cG$ of $A\rtimes_{\alpha,red}\cG$ is essential.

\subsubsection{Invariant open sets}
 Let $U$ be an open subset of $\cG$, which is saturated for $\cG$ (\ie for all $\gamma\in \cG$, we have $s(\gamma)\in U\iff r(\gamma)\in U$). Put $F=X\setminus U$. Define the subgroupoids $\cG_U=s^{-1}(U)=r^{-1}(U)$ and $\cG_F=s^{-1}(F)=r^{-1}(F)$. The action $\alpha $ of $\cG$ on $A$ gives actions $\alpha _U$ of $\cG_U$ on $A_U$ and $\alpha _F$ of $\cG_F$ on $A_F$. We may note that $A_U\rtimes _{\alpha_U}\cG_U=A_U\rtimes _{\alpha}\cG$ and $A_F\rtimes _{\alpha_v}\cG_F=A_F\rtimes _{\alpha}\cG$. Let us quote some results that we will use:
 
 \begin{enumerate}
\item We have an exact sequence of full crossed products: $$0\to A_U\rtimes _{\alpha_U}\cG_U\to  A\rtimes_\alpha \cG \to A_F\rtimes_{\alpha_F} \cG_F\to 0.$$

\item  If $\cG_F$ is amenable, the same is true for the reduced crossed products - exactness at the middle terms follows from the diagram 
$$\xymatrix{
0\ar[r]&A_U\rtimes _{\alpha_U}\cG_U\ar[r]\ar[d]  &A\rtimes_\alpha \cG\ar[r]\ar[d]& A_F\rtimes_{\alpha_F} \cG_F\to 0\ar[r]\ar[d]^{\simeq}& 0\\
0\ar[r]& A_U\rtimes_{\alpha_U,red} \cG_U\ar[r]&A\rtimes_{\alpha,red} \cG\ar[r] &A_F\rtimes_{\alpha_F,red} \cG_F\ar[r]& 0 
}
$$ 
 where the first line is exact and the vertical arrows are onto, the last one being an isomorphism.
 
\item If $A_U$ is an essential ideal in $A$, then $A_U\rtimes_{\alpha_U,red} \cG_U$ is an essential ideal in $A\rtimes_{\alpha,red} \cG$. 

\item It follows from Rem. \ref{essential} that, if $\cG_F$ is amenable and $A_U$ is an essential ideal in $A$, then $A_U\rtimes_{\alpha_U} \cG_U$ is an essential ideal in $A\rtimes_{\alpha} \cG$. 
\end{enumerate}

\section{Baaj's pseudodifferential extension}

In this section, we briefly review Baaj's construction of the pseudodifferential extension of a crossed product by a Lie group $G$. We note that the dual action extends to the pseudodifferential extension (and is trivial at the symbol level) and discuss the corresponding crossed product. Although this is not necessary in our framework, we will not assume $G$ to be abelian, so that this dual action is a coaction of $G$, since this doesn't really add any difficulty. We then establish an isomorphism between the crossed product of the algebra of the pseudodifferential operators by the dual action and a natural pseudodifferential extension. Finally, we examine the case where the Lie group is $\R$ - which is the relevant case for our results of section 5.

\subsection{Baaj's pseudodifferential calculus for an action of a Lie group}

Let us begin by recalling the extension of pseudodifferential operators associated with a continuous  action $\alpha$ by automorphisms of a Lie group $G$ on a $C^*$-algebra $A$ (\cite{BaajOPD1, BaajOPD2}, the results of Baaj concern the case $G=\R^n$ - but immediately generalize to the general case of a Lie group). 

Recall first that the order $0$ pseudodifferential operators on a Lie group $G$ give rise to an exact sequence $$0\to C^*(G) \longrightarrow \Psi^*(G)\overset{\sigma }\longrightarrow C(S^*\gag)\to 0$$ where {$C^*(G)$ is the (full) group $C^*$-algebra of $G$ and} $S^*\gag$ denotes  the (compact) space of half lines in the dual space $\gag^*$ of the Lie algebra $\gag$.

Now, the algebras $A$ and $C^*(G)$ sit in the multiplier algebra of $A\rtimes _\alpha G$ in a non degenerate way, and the elements $ax$ with $a\in A$ and $x\in C^*(G)$ span a dense subspace of $A\rtimes _\alpha G$. This holds for the full group algebra and crossed product, as well as for the reduced group algebra and crossed product. Note however that, at the level of full $C^*$-algebras, the morphism $C^*(G)\to \cM(A\rtimes  _\alpha G)$ needs not be injective in general - it is easily seen to be injective at the level of reduced $C^*$-algebras.  We will somewhat abusively  identify $C^*(G)$ and $A$ with their images in the multiplier algebra $\cM(A\rtimes  _\alpha G)$. 

In what follows, since we will consider the crossed product by the dual action, we will mainly use the reduced crossed product. Note also that we will mainly use Baaj's construction in the case where $G$ is $\R$ which is amenable and there is no distinction between the full and the reduced case. In particular the morphism $C^*(G)\to \cM(A\rtimes  _\alpha G)$ is injective in that case (if $A\ne \{0\}$).

The nondegenerate morphism $C^*(G)\to \cM(A\rtimes_{\alpha} G)$ extends to the multiplier algebra of $C^*(G)$ and in particular to the subalgebra $\Psi^*(G)$ of order $0$ pseudodifferential operators of $G$. We still identify (abusively) the elements of $\Psi^*(G)$ with their images in $\cM(A\rtimes_{\alpha} G)$. Recall that we have:

\begin{proposition}\cite[{section 4}]{BaajOPD1}
\begin{enumerate}
\item For every $P\in \Psi^*(G)$ and $a\in A$, the commutator $[P,a]$ belongs to $A\rtimes _\alpha G$.
\item The closure of the linear span of products of the form $Pa$ with $P\in \Psi^*(G)$ and $a\in A$ is a $C^*$-subalgebra $\Psi^*(A,\alpha,G)\subset \cM(A\rtimes_{\alpha} G)$ and we have an exact sequence: $$0\to A\rtimes _\alpha G \longrightarrow \Psi^*(A,\alpha,G)\overset{\sigma_\alpha}\longrightarrow C(S^*\gag) \otimes A\to 0.\eqno(1)$$
\end{enumerate}
\end{proposition}

Let us briefly discuss some naturality properties of this construction:

\begin{proposition}\label{etoile}
Let $(A,G,\alpha)$ and $(B,G,\beta)$ be $C^*$-dynamical systems and $\gamma:A\to \cM(B)$ a $G$-equivariant morphism
\begin{enumerate}
\item We obtain a morphism $\widehat \gamma:\Psi ^*(A,\alpha,G)\to \cM(\Psi^*(B,\beta,G))$ and a commutative diagram $$\xymatrix{
\Psi ^*(A,\alpha,G)\ar[r]^{\sigma_\alpha}\ar[d]^{\widehat {\gamma}} &C(S^*\gag) \otimes A\ar[d]^{\id\otimes \gamma}\\
\cM(\Psi^*(B,\beta,G))\ar[r]^{\widetilde{\sigma_\beta}} &\cM(C(S^*\gag) \otimes B) .
}
$$ both for the full and the reduced versions - where we denoted by $\widetilde{\sigma_\beta}$ the extension of $\sigma_\beta$ to the multipliers.
\item If $\gamma(A)\subset B$ then $\widehat \gamma(\Psi ^*(A,\alpha,G)) \subset \Psi^*(B,\beta,G)$. 
Moreover, if $\gamma:A\to B$ is an isomorphism, then $\widehat \gamma:\Psi ^*(A,\alpha,G)\to \Psi^*(B,\beta,G)$ is an isomorphism.
\item If $\gamma$ is injective then so is the reduced version of $\widehat \gamma$.
\end{enumerate}
\begin{proof}
\begin{enumerate}
\item By construction the inclusion of $B$ in $\Psi^*(B,\beta,G)$ is a nondegenerate morphism (\ie $B\Psi^*(B,\beta,G)=\Psi^*(B,\beta,G)$). It therefore extends to a morphism $\cM(B)\to \cM(\Psi^*(B,\beta,G))$. In this way, we find a representation $\widehat \gamma:A\to \cM(\Psi^*(B,\beta,G))$. Now the images of $A$ and $G$ in $\cM(B\rtimes _{\beta}G)\supset  
\cM(\Psi^*(B,\beta,G))$ form a covariant representation so that we get a morphism $A\rtimes _{\alpha}G\to \cM(B\rtimes _{\beta}G)$ (both for the reduced and full versions of the crossed products). The image of this morphism is spanned by elements $a.h$ with $a\in A$ and $h\in C^*(G)$; it therefore sits in $\cM(\Psi^*(B,\beta,G))$. Finally, upon replacing $A$ by the algebra obtained by adjoining a unit, we may assume that $\gamma$ is non degenerate. It follows that $\widehat \gamma :A\rtimes _\alpha G\to \cM(B\rtimes_\beta G)$ is non degenerate and therefore uniquely extends to the multiplier algebra. We thus get a morphism $\widehat \gamma :\Psi^*(A,\alpha,G)\to \cM(B\rtimes _{\beta}G)$. The image of $a.P$ is $\widehat\gamma(a).P$ (for $a\in A$ and $P\in \Psi^*(G)$) and therefore $\widehat \gamma (\Psi^*(A,\alpha,G))\subset \cM(\Psi^*(B,\beta,G))$.

\item This is obvious.
\item If $\gamma$ is one to one, then the reduced version $\gamma_{red}:A\rtimes _{\alpha,red}G\to \cM(B\rtimes_{\beta,red} G)$ is injective. Therefore   $\ker \widehat \gamma_{red}\cap A\rtimes _{\alpha, red}G=\{0\}$ whence $\ker \widehat \gamma_{red}=\{0\}$ since $A\rtimes _{\alpha,red}G$ is an essential ideal in $\Psi_{red}^*(A,\alpha,G)$ - see prop. \ref{orme}. \qedhere
\end{enumerate}
\end{proof}
\end{proposition}

\subsection{The dual action}

We now restrict to the reduced group algebras and crossed products.

The coproduct of $C_{red}^*(G)$ is a non degenerate morphism $\delta:C_{red}^*(G)\to \cM(C_{red}^*(G)\otimes C_{red}^*(G))$. It therefore extends to a morphism $\tilde \delta:\cM(C_{red}^*(G))\to \cM(C_{red}^*(G)\otimes C_{red}^*(G))$.

\begin{proposition}
The restriction of $\tilde \delta $ to $\Psi^*(G)$, is a coaction: for $P\in \Psi_{red}^*(G)$ and $f\in C_{red}^*(G)$, we have $\tilde\delta(P)(1\otimes f)\in \Psi_{red}^*(G)\otimes C_{red}^*(G)$ and the span of such products is dense in $\Psi_{red}^*(G)\otimes C_{red}^*(G)$. Moreover, for $P\in  \Psi_{red}^*(G)$ and $f\in C_{red}^*(G)$, we have $(\tilde\delta(P)-P\otimes 1)(1\otimes f)\in C^*(G\times G)$.
\begin{proof} Let $(X_i)_{1\le i\le d}$ be an (orthonormal) basis of $\gag$ and let $\Delta =-\sum_i X_i^2$ be the associated (positive) laplacian, seen as an unbounded (elliptic, positive) multiplier of $C_{red}^*(G)$.

The non degenerate morphism  $\delta $ has an extension $\check\delta$ to unbounded multipliers:
for $1\le i\le d$, set $p_i=X_i(1+\Delta)^{-1/2}\in \Psi_{red}^*(G)$.\\  We let now $C^*_{red}(G\times G)$ act faithfully on $L^2(G\times G)$. The following equalities hold on the infinite domain of the laplacian of the group $G\times G$, which is a dense subspace of $L^2(G\times G)$.

We have $\check\delta (X_i)=X_i\otimes 1+1\otimes X_i$. It follows that $\check\delta (\Delta)=\Delta\otimes 1+1\otimes \Delta-2\sum_iX_i\otimes X_i $. For $f\in C_c^\infty (G)$ {(acting as a convolution operator)}, we may then write: $$(1\otimes f)(\tilde \delta (p_i)-p_i\otimes 1)=(1\otimes fX_i)\delta((1+\Delta)^{-1/2})+(X_i\otimes f)(\delta((1+\Delta)^{-1/2})-(1+\Delta)^{-1/2}\otimes 1).$$
Now $fX_i$ and $(1+\Delta)^{-1/2}$ extend to elements of $C_{red}^*(G)$ therefore $C_i=(1\otimes fX_i)\delta((1+\Delta)^{-1/2})$ extends as well to an element of $C_{red}^*(G\times G)$. We write $(1+\Delta)^{-1/2}$ as an integral (\cf \cite{BJcras}):  $$(1+\Delta)^{-1/2} =\frac{2}{\pi}\int_0^{+\infty}(1+\Delta+\lambda^2)^{-1} d\lambda .$$ Write also 
$$(1+\Delta+\lambda^2)^{-1}\otimes 1-\delta(1+\Delta+\lambda^2)^{-1}=((1+\Delta+\lambda^2)^{-1}\otimes 1)(1\otimes \Delta+2\sum_jX_j\otimes X_j) \delta(1+\Delta+\lambda^2)^{-1}$$
Putting $D_i=(X_i\otimes f)\Big((1+\Delta)^{-1/2}\otimes 1-\delta((1+\Delta)^{-1/2})\Big)$, we find
$$\begin{array}{ccl}
D_i&=&\frac{2}{\pi}(X_i\otimes f)\int_0^{+\infty}((1+\Delta+\lambda^2)^{-1}\otimes 1)-\delta(1+\Delta+\lambda^2)^{-1}d\lambda \\
&=&\frac{2}{\pi}  \int_0^{+\infty}(X_i(1+\Delta+\lambda^2)^{-1}\otimes  f\Delta)\delta(1+\Delta+\lambda^2)^{-1}d\lambda\\&&-\frac{4}{\pi}\sum_j\int_0^{+\infty}(X_i(1+\Delta+\lambda^2)^{-1}X_j \otimes fX_j)\delta(1+\Delta+\lambda^2)^{-1}d\lambda
\end{array}$$
Now all the terms appearing are bounded operators:\begin{itemize}
\item $X_i(1+\Delta+\lambda^2)^{-1}$ is pseudodifferential of order $-1$ and therefore $X_i(1+\Delta+\lambda^2)^{-1}\in C_{red}^*(G)$;
\item $f\Delta$ and $fX_j$ are  smoothing therefore in $C_{red}^*(G)$;
\item $(1\otimes fX_j)\delta(1+\Delta+\lambda^2)^{-1}\in C_{red}^*(G)\otimes C_{red}^*(G)$.
\end{itemize}
It follows that the integrand extends to an element of $C_{red}^*(G)\otimes C_{red}^*(G)$.

Furthermore, $X_k(1+\Delta+\lambda^2)^{-1/2}=X_k(1+\Delta)^{-1/2}h_\lambda(\Delta)$ where $\|h_\lambda\|_\infty\le 1$, whence $\|X_i(1+\Delta+\lambda^2)^{-1}\|$ and $\|X_i(1+\Delta+\lambda^2)^{-1}X_j\|$ are bounded independently of $\lambda$. Hence, this integral is norm convergent and $D_i$ extends to an element $\bar D_i$ of $C_{red}^*(G)\otimes C_{red}^*(G)$.\\ Thus, we have proved that $(1\otimes f)(\tilde \delta (p_i)-p_i\otimes 1)=C_i+\bar D_i$ belongs to $C_{red}^*(G)\otimes C_{red}^*(G)$.

The set $\cA$ of $P\in \Psi_{red}^*(G)$ such that $(1\otimes C_{red}^*(G))(\tilde \delta(P)-P\otimes 1)\subset C_{red}^*(G)\otimes C_{red}^*(G)$ and $(1\otimes C_{red}^*(G))(\tilde \delta(P^*)-P^*\otimes 1)\subset C_{red}^*(G)\otimes C_{red}^*(G)$ is a closed $*$-subalgebra of $\Psi_{red}^*(G)$; it contains $C_{red}^*(G)$. As $p_i+p_i^*\in C_{red}^*(G)$, it follows by the above calculation that $p_i\in \cA$. 

Since the symbols of the $p_i$'s generate a dense subalgebra of the symbol algebra $C(S^*\gag)$ we conclude that $\cA= \Psi_{red}^*(G)$.

Finally,  the closed vector span of $(1\otimes f)\tilde \delta(P)$ contains the closed vector span of $(1\otimes f)\delta (h)$ (with $f,\ h\in C^*(G)$) hence, $C_{red}^*(G)\otimes C_{red}^*(G)$. Therefore  $(1\otimes f)\tilde \delta(P)-P\otimes f$ is in this span: the same holds for $P\otimes f$.
\end{proof}
\end{proposition}

\subsection{Isomorphisms}\label{sect3.4}

Let $\alpha$ be an action of a Lie group $G$ on a $C^*$-algebra $A$. Denote by $\hat \alpha$ the dual action on the reduced crossed product $A\rtimes_{\alpha,red}G$ as well as its extension to $\Psi_{red}^*(A,\alpha,G)$ discussed above. {Recall that in the context on non abelian groups, $B \rtimes  {\widehat G}$ is just a notation for the crossed product by a dual action, - it is a $C^*$-algebra generated by products $bf$ with $b\in B$ and $f\in C_0(G)$ subject to the equivariance condition.}

The Takesaki-Takai duality (\cite{Takai}) for non abelian groups, (see \cite{Landstad1, Landstad2}), is an isomorphism $(A\rtimes_{\alpha,red}G)\rtimes _{\hat \alpha}\widehat G\simeq A\otimes \cK$ which is based on the following facts: 
\begin{enumerate}
\item There are natural morphisms of the $C^*$-algebras $A$ and $C_0(G)$ to the multiplier algebra $\cM((A\rtimes_{\alpha,red}G)\rtimes _{\hat \alpha}\widehat G)$, as well as a (strictly continuous) morphism of the group $G$ to the unitary group of this multiplier algebra, yielding a morphism of $C^*_r(G)$ to $\cM((A\rtimes_{\alpha,red}G)\rtimes _{\hat \alpha}\widehat G)$.  

The double crossed product $(A\rtimes_{\alpha,red}G)\rtimes _{\hat \alpha}\widehat G$ is generated by the products $f.a.h$ with $a\in A$, $h\in C^*_r(G)$ and $f\in C_0(G)$ (sitting in the multiplier algebra of $(A\rtimes_{\alpha,red}G)\rtimes _{\hat \alpha}\widehat G$). Now, since the dual action is trivial on $A$, the images of $A$ and $C_0(G)$ commute so that we find in the multiplier algebra of $(A\rtimes_{\alpha,red}G)\rtimes _{\hat \alpha}\widehat G$ a copy of the $C^*$-tensor product $A\otimes C_0(G)$. The group $G$ acts on $A\otimes C_0(G)$ through the action $\alpha \otimes \lambda$ (where $\lambda$ denotes the action of $G$ on $C_0(G)$ by left translation). 

The morphisms of the $C^*$-algebra $A$ and the group $G$ (\resp of $C_0(G)$ and $G$) to $\cM((A\rtimes_{\alpha,red}G)\rtimes _{\hat \alpha}\widehat G)$ form a covariant representation of the $C^*$-dynamical system $(A,G,\alpha)$ (\resp $(C_0(G),G,\lambda)$). It follows that the morphisms of $A\otimes C_0(G)$ and $G$ in the multiplier algebra $\cM((A\rtimes_{\alpha,red}G)\rtimes _{\hat \alpha}\widehat G)$ form a covariant representation of the $C^*$-dynamical system $(A\otimes C_0(G),G,\alpha \otimes \lambda)$.

In this way, we get an isomorphism $(A\rtimes_{\alpha,red}G)\rtimes _{\hat \alpha}\widehat G\simeq (A\otimes C_0(G))\rtimes_{\alpha\otimes \lambda,red}G$.

\item Now, on $A\otimes C_0(G)$, the actions $\alpha\otimes\lambda$ and $\id\otimes \lambda$ are conjugate through the automorphism $\gamma $ of $C_0(G;A)=A\otimes C_0(G)$ given by the formula $(\gamma f)(x)=\alpha_x(f(x))$ for $f\in C_0(G;A)$ and $x\in G$. We find an isomorphism $(A\otimes C_0(G))\rtimes_{\alpha\otimes \lambda,red}G\simeq (A\otimes C_0(G))\rtimes_{\id\otimes \lambda}G$.

\item Finally $(A\otimes C_0(G))\rtimes_{\id\otimes \lambda}G\simeq A\otimes (C_0(G)\rtimes_{ \lambda}G)\simeq A\otimes \cK$.
\end{enumerate}

\begin{proposition} \label{prop2.4} The isomorphism $f:(A\rtimes_{\alpha,red} G)\rtimes _{\hat \alpha}\widehat G\overset{\sim}{\To} (A\otimes C_0(G))\rtimes _{\alpha\otimes \lambda,red}G$ extends to an isomorphism
 $\Psi_{red}^*(A,\alpha,G)\rtimes _{\hat \alpha}\widehat G\simeq \Psi_{red}^*(A\otimes C_0(G),\alpha\otimes \lambda ,G)$.
\begin{proof}
Since $A\rtimes _{\alpha,red} G$ is an essential ideal in $\Psi_{red}^*(A,\alpha,G)$ (see \ref{orme}), the algebra $\Psi_{red}^*(A,\alpha,G)\rtimes _{\hat \alpha}\widehat G$ sits in the multiplier algebra $\cM((A\rtimes_{\alpha,red} G)\rtimes _{\hat \alpha}\widehat G).$

In the same way, the algebra $\Psi_{red}^*(A\otimes C_0(G),\alpha\otimes \lambda ,G)$ sits also in $\cM((A\otimes C_0(G))\rtimes _{\alpha\otimes \lambda,red}G)$.

Both algebras are generated by products $aPh$ where $a\in A$, $P\in \Psi_{red}^*(G)$ and $h\in C^0(G)$.

Now the inclusions of $A$ and of $C_0(G)$ in $\cM$ correspond to each other under the extension $\tilde f$ of $f$ to the multipliers. As the inclusions of $C_{red}^*(G)$ to $\cM((A\rtimes_{\alpha,red} G)\rtimes _{\hat \alpha}\widehat G)$ and $\cM((A\otimes C_0(G))\rtimes _{\alpha\otimes \lambda,red}G)$ correspond to each other under $\tilde f$, the same holds for the extension to the multipliers, and in particular for the inclusions of $\Psi_{red}^*(G)$.
\end{proof}
\end{proposition}

The actions $\alpha\otimes\lambda$ and $\id\otimes \lambda$ of $G$ on  $A\otimes C_0(G)$ are conjugate.
Using prop. \ref{etoile}, we deduce isomorphisms $\Psi_{red}^*(A,\alpha,G)\rtimes _{\hat \alpha}\widehat G\simeq A\otimes \Psi_{red}^*(C_0(G),\lambda,G)\simeq A\otimes (\Psi_{red}^*(G)\rtimes _{\hat\lambda}\widehat G).$

\begin{definition}
Let $B$ be a subalgebra of $C(S^*\gag)\otimes A $. We denote by $\Psi_{red}^*(A,\alpha,G;B)$ the \emph{$B$-valued pseudodifferential extension of $\alpha$} \ie the subalgebra $$\Psi_{red}^*(A,\alpha,G;B)=\{P\in \Psi_{red}^*(A,\alpha,G);\ \sigma(P)\in B\}$$ of $\Psi_{red}^*(A,\alpha,G)$.  
\end{definition}

In the case of the trivial action, $\Psi_{red}^*(A,\id,G;B)=\{P\in A\otimes \Psi_{red}^*(G);\ ( \sigma\otimes \id )(P)\in B\}$.

\subsection{The case of $\R$}\label{caseofR}

When $G=\R$, then  $\gag^*= \R$ which has two half lines, \ie $C(S^*\gag)=\C\oplus \C$.

Extension (1) reads therefore $$0\to A\rtimes _\alpha\R \longrightarrow \Psi^*(A,\alpha,\R)\overset{\sigma_\pm}\longrightarrow A\oplus A\to 0,$$ where $\sigma_+$ and $\sigma_-$ are morphisms from $\Psi^*(A,\alpha,\R)\to A$.

It is helpful for our discussion to identify the dual group of $\R$ with $\R^*_+$ through the pairing $\langle t|u\rangle=u^{it}$ for $u\in \R_+^*$ and $t\in \R$. Under this identification, $C^*(\R)\simeq C_0(\R_+^*)$ and $\Psi_0^*(\R)\simeq C([0,+\infty])$. The maps $\sigma_-$ and $\sigma_+$ correspond to evaluation at $0$ and $+\infty$ in the sense that $\sigma_-(Pa)=P(0)a$ and $\sigma_+(Pa)=P(+\infty)a$, where $a\in A$ and $P\in C([0,+\infty])\simeq\Psi^*(\R)$.

The algebra $A$ sits in $\cM(A\rtimes_\alpha\R)$ and we have a strictly continuous family $(u_t)_{t\in \R}$ in $\cM(A\rtimes \R)$. Then we can write $u_t=Q_\alpha^{it}$ where $Q_\alpha$ is a regular unbounded, selfadjoint, positive multiplier with dense range - \ie such that $Q_\alpha^{-1}$ is also densely defined, and therefore a regular unbounded, selfadjoint, positive multiplier. The algebra $A\rtimes \R$ is spanned by $af(Q_\alpha)$ with $f\in C_0(\R_+^*)$ and $\Psi^*(A,\alpha,\R)$ is spanned by $af(Q_\alpha)$ with $a\in A$ and $f\in C([0,+\infty])$.

\begin{definition}
Let $A$ be a $C^*$-algebra and let $\alpha=(\alpha_t)_{t\in \R}$ be a continuous action of $\R$ on $A$ by $*$-automorphisms. Let $B$ be a $C^*$-subalgebra of $A$. We set $$\Psi^*(A,\alpha,\R,0,B)=\{x\in \Psi^*(A,\alpha,\R);\ \sigma_-(x)\in B,\ \sigma_+(x)=0\}.$$
\end{definition}

The algebra $\Psi^*(A,\alpha,\R,0,B)$ is spanned by elements $af(Q_\alpha)+b(1+Q_\alpha)^{-1}$ for $a\in A,\ b\in B,\ f\in C_0(\R_+^*)=C^*(\R)$ all sitting naturally as multipliers of $A\rtimes _\alpha \R$.

\section{Pseudodifferential extension associated to an action of a smooth groupoid}

In this section, we recall a few facts on smooth groupoids: the pseudodiffelential calculus, the adiabatic groupoid $\cG$ of a smooth groupoid $\cG$ \cite{MonthPie, NWX}, its ideal $J(\cG)$ (\cite[{section 4.1}]{DS1}), the action of $\R_+^*$. We then extend all these to the case of an action of $\cG$ on a $C^*$-algebra $A$. 

{Recall that $\gA \cG$ denotes the total space of the normal bundle of the inclusion of $\cG^{(0)}\subset \cG$, $\gA^*\cG$ the total space of its dual bundle, and $S^*\gA \cG$ the associated sphere bundle, \ie the set of half lines in $\gA^*\cG$.}

\subsection{The extension of pseudodifferential operators}\label{grade}

On every Lie groupoid $\cG$, there is a (longitudinal) pseudodifferential calculus. For every $m\in \R$ (and even for $m\in \C$ - \cite[{section 3}]{chief}) we have a space $\cP_m(\cG)$ of classical pseudodifferential operators of order $m$ (with polyhomogeneous symbol $\sigma\sim \sum_{k= 0}^{+\infty}a_{m-k}$ where $a_{m-k}$ is homogeneous of order $m-k$) and a symbol map which is a linear map $\sigma_m$ from $\cP_m(\cG)$ to homogeneous functions of order $m$ defined on  $\gA^*\cG$ (outside the zero section) -  with kernel $\cP_{m-1}(\cG)$. 

The smooth functions of $M=\cG^{(0)}$ define elements of $\cP_0(\cG)$; the sections of the algebroid define elements of $\cP_1(\cG)$. The algebra generated by these is the algebra of differential operators. Given a positive definite quadratic form $q$ on the bundle $\gA^*\cG$, we may find a (positive) laplacian $\Delta_\cG\in \cP_2(\cG)$ which is a positive and whose principal symbol is $q$. 

 At the level of $C^*$-algebras we obtain an extension $\Psi^*(\cG)$ of $C^*(\cG)$ and an exact sequence of order $0$ pseudodifferential operators $$0\to C^*(\cG)\longrightarrow\Psi^*(\cG)\overset{\sigma_0}\longrightarrow C(S^*\gA \cG)\to 0.$$  Recall (\cf \cite{ConnesLNM, MonthPie, NWX}) that $\Psi^*(\cG)$ is the closure of the algebra $\cP_0(\cG)$ of order zero pseudodifferential operators on $\cG$ in the multiplier algebra of $C^*(\cG)$ and $\sigma_0$ is the (extension by continuity of the) principal symbol map.

\subsection{The adiabatic groupoid and the ideal $J(\cG)$}

Let $\cG$ be a Lie groupoid. We denote by $M=\cG^{(0)}$ its set of objects. The associated \emph{adiabatic groupoid} $\cG_{ad}$  is obtained by applying the ``deformation to the normal cone'' construction to the inclusion $M\to \cG$ of the unit space of $\cG$ into $\cG$. This construction was introduced by Connes in the case of a pair groupoid $\cG=M\times M$ (\cite[{section II.5}]{ConnesNCG}), and generalized in   \cite{MonthPie, NWX}.

As a set, and as a groupoid, $\cG_{ad}=\gA \cG\times \{0\}\cup \cG\times \R_+^*$ where $\gA \cG$ is (the total space of) the Lie algebroid of $\cG$, \ie the normal bundle of the inclusion in $\cG$ of the space of objects $M$ of $\cG$; its groupoid structure is given by addition of vectors - source and range coincide and are just the bundle map $\gA \cG\to M$. These sets are glued using an exponential map $\gA \cG\to \cG$ (see \cite{MonthPie, CR, DS1} for further details).

The $C^*$-algebra of the adiabatic groupoid of $\cG$ sits in an exact sequence $$0\to C^*(\cG)\otimes C_0(\R_+^*)\longrightarrow C^*(\cG_{ad})\overset{ev_0}\longrightarrow C_0(\gA^*\cG)\to0,$$
where $\gA^*\cG$ denotes the total space of the dual bundle to the Lie algebroid $\gA \cG$ of $\cG$. Consider the morphism $\epsilon:C_0(\gA^*\cG)\to C(M)$ which associates to a function on $\gA^*\cG$ its value on the $0$-section $M$ of the bundle $\gA^*\cG$ - \ie the trivial representation of the group $\gA_x \cG$. We denote by $J(\cG)$ the kernel of $\epsilon\circ ev_0$, which is an ideal of $C^*(\cG_{ad})$. 
We therefore have an exact sequence: $$0\to J(\cG)\to C^*(\cG_{ad})\to C(M)\to 0.$$

\begin{remark}\label{Jessentiel}
It follows from \cite[Corollary 2.4]{KhoshSk}, since $M\times \R_+^*$ is dense in $M\times \R_+$ that the ideal $C_0(\R_+^*)\otimes C_{red}^*(\cG)$ is essential in $C^*_{red}(\cG_{ad})$. 

Thanks to remark \ref{essential} we deduce that $C_0(\R_+^*)\otimes C^*(\cG)$ is also an essential ideal in $C^*(\cG_{ad})$. 

As it contains $C_0(\R_+^*)\otimes C^*(\cG)$, the ideal $J(\cG)$ is essential in $C^*(\cG_{ad})$ both for the reduced and the full $C^*$-norm.

Note also that the subset $\gA^*\cG\setminus M$ is dense in $\gA^*\cG$ (unless the groupoid $\cG$ is $r$-discrete in the sense of \cite[{def. 2.6, p. 18}]{Ren} - \ie the dimension of the algebroid is $0$), and therefore $\ker \epsilon$ is essential in $C_0(\gA^*\cG)$. In this way we have another proof that $J(\cG)$ is essential in $C^*(\cG_{ad})$.
\end{remark}

We denote by $\tau$ the action of the group $\R_+^*$ by groupoid automorphisms on $\cG_{ad}$. This action is given by $\tau_t (\gamma,u) =(\gamma,tu)$ for $\gamma\in \cG$ and $t,u\in \R_+^*$ $\tau_t(x,U,0)=(x,t^{-1} U,0)$ for $ (x,U)\in \gA \cG\ \ (x\in M)$.

We therefore get an action still denoted by $\tau$ of $\R_+^*$ on $C^*(\cG_{ad})$. Note that $J(\cG)$ is invariant under this action and that the quotient action of $\R_+^*$ on $C^*(\cG_{ad})/J(\cG)=C(M)$ is trivial.

We will also use from \cite[section 3.1]{DS1} the dense subspaces $\cS(\cG_{ad})$ of $C^*(\cG_{ad}) $ and $\cJ(\cG)$ of $J(\cG)$ consisting of smooth functions with Schwartz decay properties. Recall (\cite[Theorem 3.7]{DS1}) that for $f\in \cJ(G)$ and $m\in \R$, the operator $\displaystyle\int_0^{+\infty}\! f_t\, \displaystyle\frac {dt}{t^{m+1}}$ is an order $m$ pseudodifferential operator of the groupoid $G$ \ie an element of $\cP_{m}(G)$; its principal symbol $\sigma$ is given by $\sigma(x,\xi)=\displaystyle\int_0^{+\infty} \hat f(x,t\xi,0)\displaystyle\frac {dt}{t^{m+1}}\cdot$

\subsection{Pseudodifferential extension of smooth groupoid actions}

We now extend Baaj's construction of the pseudodifferential extension to the case of an action $\alpha$ of a smooth groupoid $\cG$ on a $C^*$-algebra $A$ -  in the sense of \cite{PYLG, PYLG2} - see section \ref{essonfrerK1}. 

\subsubsection{Smooth elements}

Let $\cG$ be a smooth  groupoid with base $M$ acting on a $C_0(M)$ algebra $A$. We denote by $\alpha:s^*A\to r^*A$ this action.

We may define elements of $A$ which are smooth along the action in the following way: \begin{itemize}
\item Let   $W$ be an open subset in $ \cG$ diffeomorphic to $U\times V$ where $U\subset M$ is open and $V$ is an open ball in $\R^k$, and such that $r(u,v)=u$. Then the $C_0(W)$ algebra $(r^*A)_W$ is isomorphic to $C_0(V;A_U)$; an element $a\in r^*A$ is said to be of class $C^{\infty,0}$  if for every such $W$ and $f\in C_c^\infty (W)$, we have $fa\in C_c^\infty(V;A_U)\subset C_0(V;A_U)\simeq A_W$.
\item An element $a\in A$ is said to be smooth for the action of $\cG$ if for all $f\in C_c^\infty (\cG)$, the element $\alpha( f.(a\circ s))$ of $r^*A$ is of class $C^{\infty,0}$. Here $f.(a\circ s)$ is the class of $a\otimes f$ in $s^*A$ - \ie the restriction of $a\otimes f$ to the graph of $s$. In other words, we have $$\Big(\alpha(f.(a\circ s))\Big)_\gamma=f(\gamma)\alpha_\gamma(a_{s(\gamma)}).$$
\end{itemize}

The smooth elements form a dense sub-algebra $A^\infty$ of $A$. Indeed, if $a\in A$ and $f\in C_c^{\infty}(\cG)$, the element $f\ast a$ given by $(f\ast a)_x=\int_{G_x}f(\gamma)\alpha_\gamma a_{s(\gamma)}d\nu ^x(\gamma)$ is easily seen to be smooth. Take then a sequence $f_n$ with $f_n\in C_c^\infty(\cG)$ positive with support tending to $M$ and such that $\nu _x(f_n)=1$: we have $f_n\ast a\to a$.

\subsubsection{Crossed product by the adiabatic groupoid}\label{adgpdaction}

Let $\cG$ be a smooth  groupoid with base $M$ acting on a $C_0(M)$ algebra $A$. Consider the morphism $\cG_{ad}\to \cG\times \R_+$ which is the identity on $\cG\times \R_+^*$ and satisfies $(x,\xi,0)\mapsto (x,0)$ for $x\in M=\cG^{(0)} \subset \cG$ and $\xi \in \gag_x$. Using this morphism, the adiabatic groupoid $\cG_{ad}$ acts on the $C_0(\R_+\times M)$-algebra  $C_0(\R_+)\otimes A$: we have $A_{x,t}=A_x$ (for $t\in \R_+ $ and $x\in M$) and, for $t\in \R_+^*$, $\gamma\in \cG$ and $b\in A_{s(\gamma)}$, we have $\alpha_{\gamma,t}(b)=\alpha_\gamma(b)$; for $x\in M$, $\xi\in \gag_x$ and $b\in A_x$, we have $\alpha_{x,\xi,0}(b)=b$.

We have an exact sequence $$0\to (A\rtimes _\alpha \cG)\otimes C_0(\R_+^*)\to (A\otimes C_0(\R_+))\rtimes_\alpha \cG_{ad}\to A\otimes _{C_0(M)}C_0(\gA^*\cG)\to 0.$$
As the groupoid $\gA \cG$ is amenable, the same exact sequence holds with reduced crossed products.

Note also that the action $\tau$ of $\R_+^*$ extends on $(A\otimes C_0(\R_+))\rtimes_\alpha \cG_{ad}$: it acts naturally on $(A\otimes C_0(\R_+))=C_0(\R_+;A)$ by $(\tau_t(a))(u)=a(t^{-1}u)$.

We will also use the ideal $J(\cG,A)\subset (A\otimes C_0(\R_+))\rtimes_\alpha \cG_{ad}$ which is the kernel of the morphism $(A\otimes C_0(\R_+))\rtimes_\alpha \cG_{ad}\to A$ obtained as the composition $$(A\otimes C_0(\R_+))\rtimes_\alpha \cG_{ad}\to A\otimes _{C_0(M)}C_0(\gA^*\cG)\to A\otimes _{C_0(M)}C_0(M)=A.$$  It is the closed vector span of elements $f.a$ with $f\in J(\cG)$ and $a\in A$. It is an essential ideal in $A\otimes C_0(\R_+))\rtimes_\alpha \cG_{ad}$ (see remark \ref{Jessentiel}).

\begin{lemma}\label{ery}
If $a\in A$ is smooth for the $\cG$ action and $f\in \cS_c(\cG_{ad})$ (\cf \cite[{section} 3.1]{DS1}), then $\|[f_t,a]\|_{A\rtimes _\alpha \cG}=O(t)$.
\begin{proof}
Note that $f.a, a.f$ are in $(A\otimes C_0(\R_+))\rtimes_\alpha \cG_{ad}$ and since they are equal in $A\otimes _{C_0(M)}C_0(\gA^*\cG)$, we find that $\|[f_t,a]\|_{A\rtimes _\alpha \cG}\to 0$.

Let  $\theta :V'\to V$ be an ``exponential map'' which is a diffeomorphism of a (relatively compact) neighborhood $V'$ of the $0$ section $M$ in $\gA \cG$ onto a tubular neighborhood $V$ of $M$ in $\cG$. We assume that $r(\theta(x,U))=x$ for $x\in M$ and $U\in \gA_x \cG$. Let $W'=\{(x,U,t)\in \gA\cG\times\R_+;\ (x,tU)\in V'\}$ and $W$ be the open subset $W=\gA\cG\times\{0\}\cup V\times \R_+^*$  of $\cG_{ad}$; finally let $\Theta:W'\times \R_+\to W$ be the diffeomorphism defined by $\Theta(x,U,0)=(x,U,0)$ and $\Theta(x,U,t)=(\theta(x,tU),t)$.

If $f\in \cS_c (\R_+^*\times \cG)$, then we have $\|[f_t,a]\|_{A\rtimes _\alpha \cG}=O(t^n)$ for all $n$.

We may therefore assume that $f$ is of the form $g\circ \Theta  $ where $g\in \cS_c(W')$; then $[f_t,a]$ is the image in $A\rtimes _{\alpha}\cG$ of the function $b_t\in r^*A$, where $(b_t)_\gamma=f_t(\gamma)\Big(a_{r(\gamma)}-\alpha_\gamma\big(a_{s(\gamma)}\big)\Big)$.

Note that there is a well defined element $c\in (r\circ\Theta)^*(A\otimes C_0(\R_+))$ given by $c_{(x,U,t)}= g(x,U,t)\frac1t\Big(a_{x}-\alpha_{\theta(x,tU)}\big(a_{s(\theta(x,tU))}\big)\Big)$ for $t\ne 0$ and $-c_{(x,U,0)}$ is the derivative at $0$ of $t\mapsto \alpha_{\theta(x,tU)}\big(a_{s(\theta(x,tU))}\big)$, and $f.(c\circ \Theta^{-1})$ gives an element  $d\in (A\otimes C_0(\R_+))\rtimes_\alpha \cG_{ad}$; we have $td_t=[f_t,a]$.
\end{proof}
\end{lemma}

\subsubsection{Pseudodifferential extension}

\begin{proposition}\label{orme}
 \begin{enumerate}
\item \label{petita} For $P\in \Psi^*(\cG)$ and $a\in A$ sitting in $\cM(A\rtimes _\alpha\cG)$, we have $[P,a]\in A\rtimes_\alpha \cG$.
\item The closed vector span of products $aP$ where $a\in A$ and $P\in \Psi^*(\cG)$ is a $C^*$-subalgebra $\Psi^*(A,\alpha,\cG)\subset \cM(A\rtimes _\alpha \cG)$.
\item We have an exact sequence $$0\to A\rtimes _\alpha \cG \longrightarrow \Psi^*(A,\alpha,\cG)\overset{\sigma_\alpha}\longrightarrow A\otimes _{C_0(M)}C(S^* \gA \cG)\to 0.$$
\end{enumerate}
\begin{proof}
 \begin{enumerate}
\item We can assume $P$ is in a dense subalgebra of $\Psi^*(\cG)$ and $a$ smooth. Whence, by \cite[Theorem 3.7]{DS1}, we may choose $P=\int_0^{+\infty} f_t\frac {dt}t$ where $f=(f_t)\in \cJ(\cG)$. Then, by Lemma \ref{ery}, $[P,a]$ is a norm converging integral of elements in $A\rtimes _\alpha \cG$.
\item This closed subspace contains $A\rtimes_\alpha \cG$ and its image in $\cM(A\rtimes_\alpha \cG)/(A\rtimes_\alpha \cG)$ is a $C^*$-algebra since $\Psi^*(\cG)$ and $A$ commute in this quotient.
\item Using (\ref{petita}) and the compatibility of the inclusions of $C_0(M)$ in $\Psi^*(\cG)$ and in $\cM(A)$, we find a morphism $\varpi:C(S^* \gA \cG)\otimes_{C_0(M)} A\to \cM(A\rtimes_\alpha \cG)/(A\rtimes_\alpha \cG)$ such that $\varpi(\sigma (P)\otimes a)$ is the class of $Pa$. We just have to show that $\varpi$ is injective.

Equivalently, we wish to show that $A\rtimes_\alpha \cG$ is an essential ideal in the fibered product $\widetilde \Psi^*(\cG;A)=\Psi^*(\cG;A)\times_{\varpi(C(S^* \gA \cG))}C(S^* \gA \cG)$.

We have a representation of $\widetilde \Psi(\cG,A)$ as multipliers of $J(\cG,A)$ given, for  $(T,\sigma)\in \widetilde \Psi^*(\cG;A)$,  by $((T,\sigma)f)_t=Tf_t$ for $t\ne 0$ and $\widehat{((T,\sigma)f)_0)}(x,\xi)=\sigma(x,\xi)\widehat{f_0}(x,\xi)$,  where $T\in \Psi^*(\cG,A)$ and $\sigma\in C(S^* \gA \cG)$. This representation is faithful: indeed, if $(T,\sigma)$ is in its kernel, taking its value at $0$ it follows that  $\sigma =0$; therefore $T\in A\rtimes _\alpha G$; but the representation of $A\rtimes _\alpha \cG$ in $J(\cG;A)$ is faithful since $A\rtimes_\alpha \cG \otimes C_0(\R_+^*)\subset J(\cG,A)$. 

Now as $C_0(\R_+^*)\otimes A\rtimes_\alpha \cG$ is an essential ideal in $J(\cG;A)$, it follows that the representation  $P\mapsto 1\otimes P$ of $\widetilde \Psi^*(\cG;A)$ on $C_0(\R_+^*)\otimes A\rtimes_\alpha \cG$ is faithful, whence $A\rtimes_\alpha \cG$ is essential in $\widetilde \Psi^*(\cG;A)$.
\qedhere
\end{enumerate}
\end{proof}
\end{proposition}

\section{Action of the adiabatic groupoid and pseudodifferential extension}

Let $\cG$ be a smooth groupoid acting on the $C^*$-algebra $A$. In this section we prove the main results of this paper:\begin{itemize}
\item  We construct an action of $\R$ on the associated $C^*$-algebra $\Psi^*(\cG,A)$ of pseudodifferential operators - extending a construction sketched in \cite[{Remark 4.10}]{DS1}.
\item We establish the isomorphism  $J(\cG,A)\simeq \Psi^*(\cG,A)\rtimes_\beta \R$  - which was sketched in \cite[{Remark 4.10}]{DS1} in the case where $A=C_0(M)$ and the action is trivial.
\item Finally we identify $(A\otimes C_0(\R_+))\rtimes _{\widetilde \alpha}G_{ad}$ as a pseudodifferential extension of the above crossed product.
\end{itemize}

\subsection{The unbounded multiplier $D$ of $C^*(\cG_{ad})$}\label{vedere}

We first recall the construction of an unbounded multiplier $D$ of $C^*(\cG_{ad})$ which was given in \cite[{section} 4.4]{DS1}. 

\bigskip 
Let $\cG$ be a longitudinally smooth groupoid with compact space of objects $M=\cG^{(0)}$. 

Fix a metric on $\gA \cG$ (and therefore on $\gA^* \cG$) and choose a positive invertible pseudodifferential operator $D_1$ on $\cG$ with principal symbol $\sigma_{D_1}(x,\xi)=\|\xi\|$. It is shown in \cite[{Prop. 21}]{chief} that $D_1$ is a regular multiplier of $C^*(\cG)$. 

\begin{proposition} (\cf  \cite[Prop. 4.8]{DS1})\label{afer}
Let  $\cG$ be a Lie groupoid with compact set of objects $\cG^{(0)}=M$ and $\cG_{ad}$ its adiabatic groupoid. Fix a metric on $\gA \cG$  (and therefore on $\gA^* \cG$) and choose a positive invertible pseudodifferential operator $D_1$ on $\cG$ with principal symbol $\sigma_{D_1}(x,\xi)=\|\xi\|$.  There is a unique regular unbounded multiplier $D$ of $C^*(\cG_{ad})$ satisfying:\begin{enumerate}\renewcommand\theenumi{\roman{enumi}}
\renewcommand\labelenumi{\rm ({\theenumi})}
\item the evaluation at $1$ of $D$ is $D_1$;
\item we have $\beta_u(D)=uD$ for $u\in \R^*_+$.
\end{enumerate}
Moreover, \begin{enumerate}
\item The evaluation at $0$ of $D$, $D_0$, is the unbounded multiplier $q$ of $C_0(\gA ^*\cG)=C^*(\gA \cG)$ where $q(x,\xi)=\|\xi\|$.\label{propafera}
\item The multiplier $(1+D)^{-1}$ is in fact a strictly positive element of $C^*(\cG_{ad})$.\label{propaferb}
\item For all $f\in C_0(\R_+^*)$ we have $f(D)\in J(\cG)$. Moreover, the representation $f\mapsto f(D)$ is non degenerate: if $h\in C_0(\R_+^*)$ is strictly positive in $\R_+^*$, then $f(D)$ is a strictly positive element of $J(\cG)$..\label{propaferc}
\end{enumerate}
\begin{proof} 
If $D$  satisfies (i) and (ii), then $D_u=uD_1$ for all $u>0$, and this establishes uniqueness of $D$.

Choose a finite family $(X_1,\ldots,X_m)$ of sections of $\gA \cG$ in such a way that the embedding $\xi\mapsto \langle X_i|\xi\rangle$ is an isometry from $\gA^*\cG$ to the trivial bundle. In \cite[prop. 4.8]{DS1}, we constructed  an unbounded multiplier, call it $\widetilde D$ such that $\widetilde D_1=\Big(\sum X_i^*X_i+1\Big)^{1/2}$, $\widetilde D_0=q$ and $\widetilde D_u=u\widetilde D_1$ for $u\in \R_+^*$. Now, $D_1-\widetilde D_1$ is a $0$-order operator, whence bounded. We may then define an unbounded multiplier $D$ by putting $D_u=\widetilde D_u+u(D_1-\widetilde D_1)$ and $D_0=\widetilde D_0$.

Let us prove property (\ref{propaferb}).\\
Let $c\in \R_+^*$. Since $M\times [0,c]$ is compact and  $D$ is elliptic of order $1$ (\cite[{Th. 18 and Prop. 21}]{chief}), the restriction of $(1+D)^{-1}$ to $(\cG_{ad})_{|[0,c]}$ is in $C^*(\cG_{ad})_{|[0,c]}$.
Let $m \in \R_+^*$ such that $D_1\geq m$, we have $1+D_u\geq 1+um$ and therefore $\|(1+D_u)^{-1}\| \leq (1+um)^{-1}$. It follows that $(1+D)^{-1}$ belongs to $C^*(\cG_{ad})$.

Now, $(1+D)^{-1}C^*(\cG_{ad})$ is the domain of the multiplier  $D$, whence it is dense, and $(1+D)^{-1}$ is strictly positive.

Property (\ref{propaferc})  follows from \cite[Prop. 4.8.b)]{DS1}. Note that our $D_1$ here is slightly more general than the one used there, but the same proof applies.
\end{proof}
\end{proposition}

\subsection{The Action of $\R$ on $\Psi^*(\cG,A)$}

Let $S\in \cP_{1/2}(\cG)$ be a positive elliptic pseudodifferential operator of order $1/2$ (for instance $S$ such that $\sigma_{1/2}(S)=(\sigma_2(\Delta_\cG))^{1/4}$ where $\Delta_\cG$ is a laplacian as defined in the section \ref{grade}. Denote by $\partial_{S}$ the associated derivation on $\cM(A\rtimes_\alpha \cG)$ (see appendix - facts \ref{derivations}).

\begin{lemma}\label{gueule}
Every smooth element $a\in A$ and every classical pseudodifferential $P$ on $\cG$ of order $0$ are in the domain  of the derivation $\partial_{S}$.
\begin{proof}
We may write $S=R+S_1$ where $S_1=\int_0^{+\infty}f_t\,t^{-3/2}\,dt$,  $(f_t)$ is a positive element in $\cJ(\cG)$ and $R\in \cP_{-1/2}(\cG)$. This integral means that $\dom\,S_1$ is the set of $x\in A\rtimes_\alpha \cG$ such that the integral $\int_{0}^{+\infty}f_tx\,t^{-3/2}\,dt $ converges in norm to some $y\in A\rtimes_\alpha \cG$ and then $S_1x=y$. (Indeed, by \cite[{Prop. 21}]{chief}, $S_1$ is selfadjoint regular and it is clear that $(x,y)$ is then in the graph of $S_1^*$).

Since $a$ is assume to be smooth, the integral $\int_0^{+\infty}[f_t,a]\,t^{-3/2}\,dt$ converges in norm (by Lemma \ref{ery})  to some element $b\in A\rtimes_\alpha \cG$. Then, for $x\in \dom\,S=\dom\,S_1$, the sequence $$\int_{1/n}^{+\infty}f_tax\,t^{-3/2}\,dt=\int_{1/n}^{+\infty}(af_t+[f_t,a])x\,t^{-3/2}\,dt$$ converges in norm to $aS_1x+bx$. Thus $ax\in \dom\,S_1=\dom\, S$. It follows that $a\in \dom\, \partial_S$ and $\partial_S(a)=b+[R,a]$.

If $P\in \cP_0(\cG)$, the operator $(S^2+1)^{1/2}P(S^2+1)^{-1/2}\in \Psi^*(\cG)$; it follows that $P\dom\,S\subset \dom\,S$. Moreover, $[S,P]\in \cP_{1/2}(\cG)$ and since $\sigma_{1/2}[S,P]=[\sigma_{1/2}(S),\sigma_{1/2}(P)]=0$, we find $[S,P]\in \cP_{-1/2}(\cG)\subset C^*(\cG)$.
\end{proof}
\end{lemma}

\begin{proposition}\label{ephant}
Let $D_1\in \cP_1(\cG)$ be any positive invertible pseudodifferential operator elliptic of order $1$. Then we have an action $\beta$ of $\R$ on $\Psi^*(\cG;A)$ given by $\beta_t(P)=D_1^{it}PD_1^{-it}$. This action is trivial at the symbol level.
\begin{proof}
By \cite[{Theorem 41}]{chief} there exists $S\in \cP_{1/2}$  positive elliptic of order $1/2$ and $T\in C^*(\cG)$ such that $\sqrt{D_1}=S+T$. It follows by Lemma \ref{gueule}, that with $a,P$ as above  $Pa\in \dom\,\partial_{\sqrt{D_1}}$.

Since $D_1^{-1/2}\in A\rtimes_\alpha \cG$, it follows from Lemma \ref{LemmeApp3}, that $Pa\in \dom\,\partial_{\ln D_1}$ and $[\ln D_1,aP]\in A\rtimes _\alpha \cG$. The conclusion follows from Lemma \ref{LemmeApp1}.
\end{proof}
\end{proposition}

\subsection{Isomorphism $\Psi^*(\cG,A)\rtimes \R\simeq J(\cG,A)$}\label{blonde}

In \cite[prop. 4.2.b)]{DS1}, we constructed a morphism $\phi:\Psi^*(\cG)\to \cM(J(\cG))$ such that, for $P\in\Psi^*(\cG)$ and $f=(f_u)$ in $J(\cG)$ we have  $(\phi(P)(f))_u=P\ast f_u$, for $u\ne 0$ and $(\phi(P)(f))_0=\sigma_0(P)f_0$ thanks to \cite[prop. 4.2.b)]{DS1}. Now $J(\cG)$ sits in a non degenerate way in $\cM(J(\cG,A))$. Also, by definition $A$ embeds in a compatible way in $\cM(J(\cG,A))$.

In this way, we find a morphism $\phi:\Psi^*(\cG,A)\to \cM(J(\cG,A))$ such that, for $P\in\Psi^*(\cG,A)$ and $f=(f_t)$ in $J(\cG,A)$ we have  $(\phi(P)(f))_u=P\ast f_u$, for $u\ne 0$ and $(\phi(P)(f))_0=\sigma_0(P)f_0$.

Furthermore, the operator $D$ recalled in section \ref{vedere} yields a one parameter group $(D^{it})_{t\in \R}$ in $\cM(J(\cG))$; we will still denote by $(D^{it})_{t\in \R}$ its image in $\cM(J(\cG,A))$.

As $D_u$ and $D_1$ are scalar multiples of each other, we find in this way a covariant representation of  the pair $(\Psi^*(\cG),\beta,\R)$ (prop. \ref{ephant}).

 Associated to this covariant representation  is a morphism from $\Psi^*(\cG,A)\rtimes_\beta \R$ into the multiplier algebra of $J(\cG)$, but since  the image of $C^*(\R)\subset \Psi^*(\cG)\rtimes_\beta \R$ is contained in $J(\cG)$, we get a homomorphism $\varphi :\Psi^*(\cG,A)\rtimes_\beta \R\to J(\cG,A)$. For $P\in \Psi^*(\cG,A)\subset \cM(\Psi^*(\cG,A)\rtimes_\beta \R)$ and $f\in C^*(\R)=C_0(\R_+^*)\subset \cM(\Psi^*(\cG,A)\rtimes_\beta \R)$, we have $(\varphi(Pf)) =\phi(P)f(D)$.
 
\begin{proposition}
The homomorphism $\varphi$ is an equivariant isomorphism from $(\Psi^*(\cG,A)\rtimes_\beta \R,\hat\beta)$ to $(J(\cG,A),\tau)$.
\begin{proof}
The images of  the elements of $\Psi^*(\cG,A)$ are translation invariant, \ie invariant by the extension $\overline\tau_u$ of $\tau_u$ to the multiplier algebra, and $\overline\tau_u(D^{it})=u^{it}D^{it}$. This shows that $\varphi$ is an equivariant morphism from $(\Psi^*(\cG,A)\rtimes_\beta \R,\hat\beta)$ to $(J(\cG,A),\tau)$.

Now $\beta_t$ restricts to an action of $\R$ on $C^*(\cG)$, and according to \cite[prop. 4.2.a)]{DS1} it follows that $\varphi$ extended to the multipliers defines a morphism from $C^*(\cG)\rtimes_\beta \R$ into the ideal $C_0(\R_+^*)\otimes C^*(\cG)$ of $J(\cG)$. It follows that $\varphi(A\rtimes_\alpha \cG)$  is contained in the ideal $A\rtimes_\alpha \cG\otimes C_0(\R_+^*)$ of $J(\cG,A)$. We thus have the diagram: 
$$\xymatrix{
0\ar[r]&(A\rtimes _{\alpha}\cG)\rtimes_\beta \R\ar[r]\ar[d]^{\varphi'}  &\Psi^*(\cG,A)\rtimes _\beta\R\ar[r]\ar[d]^{\varphi}& (A\otimes _{C(M)} C(S^* \gA \cG))\rtimes _\beta\R\ar[r]\ar[d]^{\varphi''}& 0\\
0\ar[r]& (A\rtimes _\alpha \cG)\otimes C_0(\R_+^*)\ar[r]&J(\cG,A)\ar[r] &A\otimes _{C(M)}C_0(\gA^*\cG\setminus M)\ar[r]& 0 
}
$$ 
As $D_1$ is an unbounded invertible multiplier of $C^*(\cG)$ and therefore of $A\rtimes _\alpha \cG$, the action $\beta$ of $\R$ on $A\rtimes _\alpha \cG$ is inner. It follows that the crossed product $(A\rtimes _\alpha \cG)\rtimes _\beta \R$ identifies with $(A\rtimes _\alpha \cG)\otimes \R_+^*$. This isomorphism is defined in the following way:  the canonical multipliers of the crossed product, \ie the generators $a\in A\rtimes _\alpha \cG$ and $\lambda_t$ for $t\in \R$ map to the functions $u\mapsto a$ and $u\mapsto u^{it}D_1^{it}$ from $\R_+^*$ to $\cM(A\rtimes _\alpha \cG)$. It follows, the image of $af$ with $a\in C^*(\cG)$ and $f\in C^*(\R)=C_0(\R_+^*)$ is $af(D)$. This isomorphism identifies thus with $\varphi'$.

The action $\beta$ is trivial on symbols; thus $(A\otimes _{C(M)}C(S^* \gA \cG))\rtimes _\beta\R$ is equal to $(A\otimes _{C(M)} C(S^* \gA \cG))\otimes C_0(\R_+^*)$, and $\varphi''(\sigma \otimes f)=\sigma f(q)$ is the isomorphism corresponding to the homeomorphism $\gA^*\cG\setminus M\simeq S^* \gA \cG\times \R_+^*$ given by $\xi\mapsto (\xi/q(\xi),q(\xi))$. The result follows. 
\end{proof}
\end{proposition}

\subsection{The crossed product by the adiabatic groupoid}

The algebra $A$ sits in $\Psi^*(\cG,A)$ as (the closure of) order $0$ differential operators. Denote by $\vartheta:A\to \Psi^*(\cG;A)$ the corresponding morphism. The element $\vartheta(a)$ as a multiplier of $A\rtimes_\alpha \cG$, is just the multiplication by $a$. 

\begin{remark}\label{mult=mult}
Using at the non degenerate morphism $\Psi^*(\cG,A)\to \cM(J(\cG;A))$ we then obtain a morphism $\hat \vartheta:A\to \cM(\Psi^*(\cG,A)\rtimes \R)$.

Also the algebra $A$ is in the multiplier algebra of $A\otimes C_0(\R_+)$ end thus we have an embedding $\tilde \vartheta :A\to \cM((A\otimes C_0(\R_+))\rtimes _{\tilde \alpha}\cG_{ad})$ - which is a subalgebra of $\cM(J(\cG;A))$ since $J(\cG;A)$ is an essential ideal in $(A\otimes C_0(\R_+))\rtimes _{\tilde \alpha}\cG_{ad}$. 
\end{remark}

We now use the notation of paragraph \ref{caseofR}. The main result of this paper is :

\begin{theorem}\label{maintheorem}
The isomorphism $\varphi:\Psi^*(\cG,A)\rtimes_\beta \R\to J(\cG,A)$ extends uniquely to an isomorphism of $\Psi^*(\Psi^*(\cG,A),\beta,\R,0,A)$ with $(A\otimes C_0(\R_+))\rtimes_{\widetilde \alpha}\cG_{ad}$. This isomorphism intertwines the actions $\beta $ and $\tau$ of $\R$.
\begin{proof}
The isomorphism $\varphi:\Psi^*(\cG,A)\rtimes _\beta\R\to J(\cG,A)$ extends to an isomorphism $\Phi$ of the multiplier algebras. Since the ideals $\Psi^*(\cG,A)\rtimes _\beta\R\subset \Psi^*(\Psi^*(\cG,A),\beta,\R)$ and $J(\cG,A)\subset (A\otimes C_0(\R_+))\rtimes_\alpha \cG_{ad}$ are essential, we just need to show that $\Phi (\Psi^*(\Psi^*(\cG,A),\beta,\R,0,A))=(A\otimes C_0(\R_+))\rtimes_{\widetilde \alpha}\cG_{ad}$.

It follows from proposition \ref{etoile}.a) that the morphism $\Phi$ coincides on $\Psi^*(\cG,A)$ with the morphism $\phi : \Psi_*(\cG,A)\rightarrow \cM(J(\cG,A))$ of  {section} \ref{blonde} and that the image of the unbounded multiplier $Q_\beta$ (see {section} \ref{caseofR}) is $D$. 

With the notation introduced in remark \ref{mult=mult}, one easily checks that $\Phi\circ \hat\vartheta =\tilde\vartheta$. 

We deduce that $\Phi\Big(\Psi^*(\Psi^*(\cG,A),\beta,\R,0,A)\Big)$ is spanned by $\varphi(\Psi^*(\cG,A)\rtimes _\beta\R)= J(\cG,A)$ and $(1+D)^{-1}\tilde\vartheta(a)$ where and $a$ over $A$.

Since $(1+D)^{-1}\in C^*(\cG_{ad})$ (prop. \ref{afer}.\ref{propaferb})), and for $a\in A$ we have $(1+D)^{-1}\tilde\vartheta(a))\in  (A\otimes C_0(\R_+))\rtimes _{\tilde \alpha}\cG_{ad})$. 

Finally $\Phi$ induces a homomorphism $\widetilde \varphi:\Psi^*(\Psi^*(\cG,A),\beta,\R,0,A))\to  (A\otimes C_0(\R_+))\rtimes_{\widetilde \alpha}\cG_{ad}$.

Moreover, since $\ev_0(D)=q$ which vanishes at the $0$ section of $\gA^* \cG$, we find that $\epsilon\circ \ev_0((1+D)^{-1})=1$, whence $\epsilon\circ \ev_0(\Phi((1+D)^{-1}\theta(a))=a$.
We thus have a commutative diagram where the sequences are exact:
$$\xymatrix{
0\ar[r]&\Psi^*(\cG,A)\rtimes_\beta \R\ar[r]\ar[d]^{\varphi}  &\Psi^*(\Psi^*(\cG,A),\beta,\R,0,A)\ar[r]\ar[d]^{\Phi}& A\ar[r]\ar[d]^{\id_A}& 0\\
0\ar[r]& J(\cG,A)\ar[r]&(A\otimes C_0(\R_+))\rtimes _{\tilde \alpha}\cG_{ad})\ar[r] &A\ar[r]& 0 
}
$$ 
Whence $\widetilde \varphi$ is an isomorphism.

By uniqueness of the extension to multipliers, we deduce that $\hat \beta_t\circ \widetilde \varphi=\widetilde \varphi\circ \tau_t$ for all $t\in \R_+^*$.
\end{proof}
\end{theorem}

Recall that the gauge adiabatic groupoid $\cG_{ga}$ is the semi-direct product $\cG_{ga}=\cG_{ad}\rtimes_\tau \R_+^*$. If $\cG$ acts on $A$, then $\cG_{ga}$ acts on $A\otimes C_0(\R_+)$.

\begin{corollary}
We have isomorphisms 
\begin{eqnarray*}
(A\otimes C_0(\R_+))\rtimes_{\alpha} \cG_{ga}&\simeq &\Psi^*(\Psi^*(\cG,A),\beta,\R,0,C(M))\rtimes_{\hat\beta}\R_+^*\\
 &\simeq &\Psi^*(\Psi^*(\cG,A)\otimes C_0(\R),\beta\otimes \lambda,\R,0,C(M)\otimes C_0(\R)).
\end{eqnarray*}
\begin{proof}
We have $ (A\otimes C_0(\R_+))\rtimes_{\alpha} \cG_{ga}= ((A\otimes C_0(\R_+))\rtimes_{\alpha} \cG_{ad})\rtimes_{\tau}\R_+^*$.
The first isomorphism is a direct consequence of theorem \ref{maintheorem}; the second one comes from prop. \ref{prop2.4}. 
\end{proof}
\end{corollary}

\begin{remark}
Let us drop the algebra $A$. The exact sequence $$0\to C^*(\cG)\otimes \cK\to C^*(\cG_{ga})\rtimes \R_+^* \to C_0(\gA^*\cG)\rtimes \R_+^*\to 0$$ defines an ``ext'' element in $KK^1(C_0(\gA^*\cG)\rtimes \R_+^*,C^*(\cG)\otimes \cK)$. Using Connes' Thom isomorphism (\cf \cite{ConnesThom, FaSkThom}), this group is isomorphic to $KK(C_0(\gA^*\cG),C^*(\cG))$. In fact, using again the Thom isomorphism, this element corresponds to the ext element  in $KK^1(C_0(\gA^*\cG),C^*(\cG)\otimes C_0(\R_+^*))$ of the exact sequence $$0\to C^*(\cG)\otimes C_0(\R_+^*)\to C^*(\cG_{ad})\to C_0(\gA^*\cG)\to 0.$$
{One easily sees (using e.g. \cite[Theorem 2.1]{MonthPie}) that} this element is the analytic index.

Let $\mu:C(M)\to \Psi^*(\cG)$ be the inclusion, and let $C_\mu$ be the corresponding mapping cone. We have an exact sequence $$0\to \Psi^*(\cG)\otimes C_0(\R_+^*)\to C_\mu\to C(M)\to 0.$$ The quotient of $C_\mu$ by the ideal $C^*(\cG)\otimes C_0(\R_+^*)$ is the cone of the inclusion $C(M)\to C(S^*\gag)$, which is naturally isomorphic to $C_0(\gA^*\cG)$. We thus find an exact sequence $$0\to C^*(\cG)\otimes C_0(\R_+^*)\to C_\mu\to C_0(\gA^*\cG)\to 0.$$
The corresponding $KK$-element can be seen again to be the analytic index element in $KK(C_0(\gA^*\cG),C^*(\cG))$.
Taking crossed product by the natural action of $\R_+^*$ on $C_\mu$ (just by rescaling), we find an exact sequence $$0\to \cK\to C_\mu\rtimes \R_+^* \to C_0(T^*M)\rtimes \R_+^*\to 0.$$

In the case of the pair groupoid, we deduce an isomorphism $C_\mu\rtimes \R \simeq C^*(\cG_{ga})$ thanks to Voiculescu's theorem (\cite[{Theorem 1.5}]{Voicu}).

It is a natural question to decide whether this isomorphism extends to the general case. On the other hand, this isomorphism is not ``natural''.  Indeed, $C_\mu$ and $C^*(\cG_{ad})$ are not isomorphic in general, whence there is no isomorphism $C_\mu\rtimes \R \simeq C^*(\cG_{ga})(=C^*(\cG_{ad})\rtimes \R_+^*)$ equivariant with respect to the dual actions. 
\end{remark}

\section{Appendix: some facts on unbounded operators}

In this appendix, we recall a few rather classical abstract facts about unbounded operators that we used in the text. These facts are presented here in a form suitable for our exposition and certainly not in their most general forms. They can be found in (or deduced directly from) \cite{BaajThese, Woro} - see also \cite{chief}.

Let $E$ be a $C^*$-module (over a $C^*$-algebra) and $L$ a regular (densely defined, unbounded) self-adjoint operator on $E$. 

\begin{facts}\label{ract}
Let us recall a few facts about unbounded functional calculus, $f\mapsto f(L)$ (\cf \cite{BaajThese, Woro}).

\begin{enumerate}
\item Put $h(t)=(i+t)^{-1}$; there exists a unique morphism $\pi_L:f\mapsto f(L)$ from $C_0(\R)$ to $\cL(E)$ such that $\pi_L(h)=(L+i\,\id_E)^{-1}$.

\item Since $h(L)$ has a dense range ($\dom\, L$), this morphism is non degenerate, it extends to a morphism $f\mapsto f(L)$ from $C_b(\R)=\cM(C_0(\R))$ to $\cL(E)$.

\item If $f\in C(\R)$, define the operator $f(L)$ whose domain is the range of $g(L)$ where $g(t)=(|f(t)|+1)^{-1}$ and such that $f(L)g(L)=(fg)(L)$.

\item If $f,g\in C(\R)$ are such that $\frac{f}{|g|+1}$ is bounded, then $\dom\,g(L)\subset \dom f(L)$.\label{ractaire}

\item \label{rigerateur}If $(f_n)$ is an increasing sequence of positive elements of $C_b(\R)$ converging simply (and therefore uniformly on compact subsets of $\R$) to a continuous function $f$, then the domain of $f(L)$ is the set of $x\in E$ such that $(f_n(L)x)$ converges (in norm) and then $f(L)x$ is the limit of this sequence.

Indeed, as $\frac{f_n+1}{f+1}=h_n$ converges to $1$ for the topology of $C_b(\R): $\begin{itemize}
\item if $x$ is in the domain of $f(L)$, it is of the form $x=(f(L)+1)^{-1}z$, and $x+f_n(L)x=h_n(L)z$ converges to $z$, therefore $f_n(L)x$ converges to $z-x$;
\item  $(f(L)+\id_E)^{-1}(f_n(L)x+x)=h_n(L)x$ converges to $x$; assume that $f_n(L)x$ converges to $y\in E$, then  $(x,x+y)$ is the limit of the sequence $\Big(h_n(L)x,(f_n(L)x+x)\Big)$ of elements of the graph of $f(L)+\id_E$; therefore $y=f(L)x$ since the graph of $f(L)$ is closed.\end{itemize}
\end{enumerate}
\end{facts}

\begin{lemma}\label{mondo}
 We have an equality $$L=\int_1^{+\infty}\Big(\frac{1}{s}-(e^L+s)^{-1}\Big)ds-\int_0^1(e^L+s)^{-1} ds$$ which means that $\dom\,L$ is the set of $x\in E$ such that the integrals $$\int_1^{+\infty}\Big(\frac{1}{s}-(e^L+s)^{-1}\Big)x \, ds\ \ \hbox{and}\ \ \int_0^1(e^L+s)^{-1} x\, ds$$ are norm convergent and $Lx$ is then the difference of these two integrals.

\begin{proof}
Put $f_n(t)=\int_1^{n}\Big(\frac{1}{s}-(e^t+s)^{-1}\Big)ds$ and $f(t)=\lim f_n(t)=\ln (e^t+1)$; put also $g_n(t)=\int_{\frac1n}^1(e^t+s)^{-1} ds$ and $g(t)=\lim g_n(t)=\ln (e^t+1)- t$. 

Then as $\frac{\ln (e^t+1)}{|t |+1}$ is bounded, $\dom\, L=\dom\, f(L)\cap \dom\, g(L)$ (by fact \ref{ract}.\ref{ractaire}).  
The conclusion follows from fact \ref{ract}.\ref{rigerateur}). 
\end{proof}
\end{lemma}

\begin{description} 
\item [Fact \ref{ract}.f).] Assume $L$ is positive with resolvent in $\cK(E)$. Then $f\mapsto f(L)$ is a morphism $\pi_L:C_0(\R_+^*)\to\cK(E)$. Note that, for $t\in \R_+^*$, we have $\pi_{tL}=\pi_L\circ \lambda_t$ where $\lambda_t$ is the automorphism of $C_0(\R_+^*)$ induced by the regular representation. Since $t\mapsto \frac{t}{t^2+1}$ is a strictly positive element of $C_0(\R_+^*)$, it follows that $\pi_L(C_0(\R_+^*))E$ is the closure of the image of $L(L^2+1)^{-1}$.
\end{description}

\begin{facts}[about derivations]
 \label{derivations}We will consider the (unbounded, skew adjoint) derivation $\partial_L$ associated with $L$: its domain is the $*$-subalgebra of the elements $a\in \cL(E)$, such that there exists $\partial _L(a)\in \cL(E)$ with $aL\subset La+\partial_L (a)$ (in other words $a\,\dom\, L\subset \dom\, L$ and $[a,L]$ defined on $\dom\, L$ extends to an operator $\partial _L(a)\in \cL(E)$). 

Put $u_t=\exp(itL)$ and define for $a\in \cL(E)$,  $\beta_t(a)=u_tau_t^*$. 

\begin{enumerate}
\item For $a\in \cL(E)$, the map $t\mapsto \beta_t$ is of class $C^1$ (for the norm topology) if and only if $a\in \dom\, \partial_L$ and, in that case $d/dt(\beta_t(a))=i\partial_L (\beta_t(a))=i\beta_t(\partial_L (a))$.

\item The closure $\overline{\dom\, \partial_L}$ of $\dom\, \partial_L$ is a $C^*$-subalgebra of $\cL(E)$ and $t\mapsto \beta_t(a)$ is a continuous action of $\R$ on it.
\end{enumerate}
\end{facts}

\begin{lemma}\label{LemmeApp1}
Let $Q$ be the norm closure of $\{a\in \dom\, \partial_L;\ \partial_L a\in \cK(E)\}$. It is a $C^*$-subalgebra of $\overline{\dom\, \partial_L}$ invariant under the action $\beta$ of $\R$. The quotient action of $\R$ on $Q/\cK(E)$ is trivial. In particular, every $C^*$-subalgebra of $Q$ containing $\cK(E)$ is invariant by $\beta$. 
\begin{proof}
Denote by $q:\cL(E)\to \cL(E)/\cK(E)$ the quotient map. If $a\in \dom\,\partial_L$ satisfies $\partial_L a\in \cK(E)$, then $t\mapsto \beta_t(a)$ is $C^1$, and the derivative of $t\mapsto q(\beta_t(a))$ is zero. All other statements are clear.
\end{proof}
\end{lemma}

\begin{lemma}\label{LemmeApp2}
Let $a\in \dom\,\partial_{e^L}\cap  \partial_{e^{-L}}$. Then $a\in \dom\,\partial_L$. 
 If the resolvent of $L$ is in $\cK(E)$, then $\partial_L(a)\in \cK(E)$.
\begin{proof}
The integral $\int_1^{+\infty}\Big[\frac{1}{s}-(e^L+s)^{-1},a\Big]ds=\int_1^{+\infty}(e^L+s)^{-1}[e^L,a](e^L+s)^{-1}\,ds$ is norm convergent (since $\|(e^L+s)^{-1}\|\le s^{-1}$), as well as 
\begin{eqnarray*}
 -\int_0^1\Big[(e^L+s)^{-1},a\Big] ds&=&\int_0^1(e^L+s)^{-1}\Big[e^L,a\Big](e^L+s)^{-1}\,ds\\
 &=&-\int_0^1e^L(e^L+s)^{-1}\Big[e^{-L},a\Big]e^L(e^L+s)^{-1}\,ds.
\end{eqnarray*}
(since $\|e^L(e^L+s)^{-1}\|\le 1$).

It follows, with the notation of Lemma \ref{mondo} that $[(f_n-g_n)(L),a]$ converges to an element $b=\int_0^{+\infty}(e^L+s)^{-1}[e^L,a](e^L+s)^{-1}\,ds$. If $x\in \dom\,L$, then $(f_n-g_n)(L)ax$ converges to $aLx+bx$; therefore $ax\in \dom \,L$ and $\partial_L(a)=b$.

Assume $L$ has compact resolvent {(\ie in $\cK(E)$)}. Put $q_s=(e^L+s)^{-1}[e^L,a](e^L+s)^{-1}$. Note that $e^Lq_s$ is bounded and, since $q_s=-(e^L+s)^{-1}e^L\Big[e^{-L},a\Big]e^L(e^L+s)^{-1}$, $e^{-L}q_s$ is also bounded. If $L$ has compact resolvant, then $(e^L+e^{-L})^{-1}\in \cK(E)$, whence $q_s\in \cK(E)$.
\end{proof}
\end{lemma}

\begin{lemma}\label{LemmeApp3}
Assume $L$ is positive. Let $a\in \cL(E)$ such that $a\,\dom\, e^L\subset \dom\, e^L$ and $e^{-L/2}[e^L,a]$ defined on $\dom\, e^L$ extends to an element of $\cL(E)$. 
Then $a\in \dom\,\partial_L$. 
 If moreover the resolvant of $L$ is in $\cK(E)$, then $\partial_L(a)\in \cK(E)$.
\begin{proof}
The integral $\int_1^{+\infty}\Big[\frac{1}{s}-(e^L+s)^{-1},a\Big]ds=\int_1^{+\infty}(e^L+s)^{-1}[e^L,a](e^L+s)^{-1}\,ds$ is norm convergent, since $$\|(e^L+s)^{-1}[e^L,a](e^L+s)^{-1}\|\le \|e^{L/2}(e^L+s)^{-1}\|\|e^{-L/2}[e^L,a]\|\|(e^L+s)^{-1}\|\le s^{-1/2}Cs^{-1}$$ for $C=\|e^{-L/2}[e^L,a]\|$. 

Of course the integral $ -\int_0^1\Big[(e^L+s)^{-1},a\Big] ds$ is also norm convergent.

It follows, with the notation of Lemma \ref{mondo} that $[(f_n-g_n)(L),a]$ converges to an element $b=\int_0^{+\infty}(e^L+s)^{-1}[e^L,a](e^L+s)^{-1}\,ds$. If $x\in \dom\,L$, then $(f_n-g_n)(L)ax$ converges to $aLx+bx$; therefore $ax\in \dom \,L$ and $\partial_L(a)=b$.

Assume $L$ has compact resolvant. Then, since $L$ is positive, $(e^L+s)^{-1}\in \cK(E)$, whence $(e^L+s)^{-1}[e^L,a](e^L+s)^{-1}\in \cK(E)$.
\end{proof}
\end{lemma}

\bibliography{BiblioTP.bib} 
\bibliographystyle{amsplain}

\end{document}